\documentclass[12pt]{article}
\usepackage{latexsym,amsfonts,amssymb}
\usepackage[dvips]{color}
\usepackage{colordvi,multicol}
\usepackage{amsmath}
\usepackage{graphicx}
\usepackage{pdflscape}
\usepackage{rotating}

\usepackage{latexsym,amsfonts,amssymb}
\usepackage[dvips]{color}
\usepackage{colordvi,multicol}
\usepackage{amsmath}
\usepackage{graphicx}
\usepackage{pdflscape}
\usepackage{rotating}
\usepackage{breqn}

\def \cal{\mathcal}

\newtheorem{thm}{Theorem}[section]

\newtheorem{lem}[thm]{Lemma}
\newtheorem{pro}[thm]{Proposition}

\newtheorem{rem}[thm]{Remark}

\newtheorem{con}[thm]{Conjecture}

\date{}
\begin{document}

\title{\bf Moment ratio inequality of bivariate Gaussian
distribution and three-dimensional Gaussian product inequality}
 \author{Oliver Russell and Wei Sun\\ \\ \\
  {\small Department of Mathematics and Statistics}\\
    {\small Concordia University, Canada}\\ \\
{\small o\_russel@live.concordia.ca,\ \ \ \ wei.sun@concordia.ca}}

\maketitle

\begin{abstract}

\noindent We prove the three-dimensional Gaussian product inequality (GPI) $E[X_1^{2}X_2^{2m_2}X_3^{2m_3}]\ge E[X_1^{2}]E[X_2^{2m_2}]E[X_3^{2m_3}]$ for any centered Gaussian random vector $(X_1,X_2,X_3)$ and $m_2,m_3\in\mathbb{N}$. We discover a novel inequality for the moment ratio $\frac{|E[ X_2^{2m_2+1}X_3^{2m_3+1}]|}{E[ X_2^{2m_2}X_3^{2m_3}]}$, which implies the 3D-GPI. The interplay between computing and hard analysis plays a crucial role in the proofs.
\end{abstract}



\section{Introduction and main results}

Correlation is a fundamental notion of science and the Gaussian distribution, which is a building block of probability and statistics, has a nearly 300-year history dating back to Abraham de Moivre. When combined, these two concepts have proven essential to many branches of mathematics and physics. However, our understanding of them is far from complete.

Inequalities involving correlations between Gaussian random variables have important connections to a number of active research areas, such as small-ball probabilities (cf. \cite{Li} and \cite{Shao}) and the zeros of random polynomials (cf. \cite{LS}). In 2014, Royen's proof of the elusive Gaussian correlation inequality conjecture \cite{Royen,LM} was a significant breakthrough that inspired renewed hope in solving another long-standing and challenging problem, the Gaussian product inequality (GPI) conjecture, which (see, for example, \cite{MNPP16}) implies the `real linear polarization constant conjecture' \cite{BST98} and is deeply linked to Kagan, Linnik and Rao's famous $U$-conjecture \cite{KLR1973}.

Most generally, the GPI conjecture (see Li and Wei \cite{LW12}) states that for any non-negative real numbers $y_j$, $j=1,\ldots,{n}$, and any ${n}$-dimensional real-valued centered Gaussian random vector $(X_1,\dots,X_{n})$,
\begin{eqnarray}\label{LW-inequ}
E \left[\prod_{j=1}^{n}|X_j|^{y_j}\right]\geq \prod_{j=1}^{n}E[|X_j|^{y_j}].
\end{eqnarray}
The difficulty of proving the GPI lies in the subtlety of the inequality and thus, so far, only partial results have been obtained. Frenkel \cite{Fr08} proved (\ref{LW-inequ}) when all exponents equal $2$ using algebraic methods. By using the Gaussian hypergeometric function as a powerful tool, Lan et al. \cite{LHS} were successful in proving the following $3$-dimensional GPI: for any $m_1,m_2\in\mathbb{N}$ and any centered Gaussian random vector $(X_1,X_2,X_3)$,
$$
E[X_1^{2m_1}X_2^{2m_2}X_3^{2m_2}]\ge E[X_1^{2m_1}]E[X_2^{2m_2}]E[X_3^{2m_2}].
$$

Recently, we gave a complete characterization of the $2$-dimensional GPI in \cite{RusSun3} and designed a computational algorithm based on sums-of-squares that can potentially be used to rigorously prove any GPI with all but one exponent fixed in \cite{RusSun2}. Furthermore, in \cite[Lemma 2.3]{RusSun} we proved a stronger version of the GPI for even-integer exponents when all correlations are non-negative (see also \cite{GenOui} for a related result, and \cite{Edel} and \cite{GenOui2} for distributional generalizations). In the same paper, we used combinatorial methods to prove the following GPIs: for any $m_3\in\mathbb{N}$ and any centered Gaussian random vector $(X_1,X_2,X_3)$,
$$
E[X_1^{2}X_2^{4}X_3^{2m_3}]\ge E[X_1^{2}]E[X_2^{4}]E[X_3^{2m_3}],
$$
and
$$
E[X_1^{2}X_2^{6}X_3^{2m_3}]\ge E[X_1^{2}]E[X_2^{6}]E[X_3^{2m_3}],
$$
and left the remainder of Theorem \ref{thmm1m2m3} below as a conjecture.

Throughout this paper, any Gaussian random variable is assumed to be real-valued and nondegenerate,
i.e., has positive variance. We will prove the following 3D-GPI.
\begin{thm}\label{thmm1m2m3}
Let $m_2,m_3 \in \mathbb{N}$. For any centered Gaussian random vector $(X_1,X_2,X_3)$,
\begin{equation}\label{GPIk}
E[X_1^{2} X_2^{2m_2} X_3^{2m_3}]\ge E[X_1^{2}] E[X_2^{2m_2}] E[X_3^{2m_3}].
\end{equation}
The equality sign holds if and only if $X_1,X_2,X_3$ are independent.
\end{thm}

In \cite{LHS}, all 3 exponents were unbounded, but 2 were equal to each other, allowing for some symmetry and  convexity to be exploited. The main contribution of this present paper lies in the fact that we consider $3$ \textit{distinct} exponents, $2$ of which are unbounded. Even this special case is extremely challenging since the subtle nature of the GPI resists even small modifications. For this reason, we used the computer (\textbf{Mathematica}) as a time-saving guide at key junctures to help us guess appropriate simplifications under which the GPI still holds. We emphasize that despite this guidance, the computer was less of a crutch and more of a shot in the arm, and all proofs contained in this paper remain rigorous and easily verifiable. We strongly believe that computing will be an essential tool used in the eventual complete proof of the GPI.

As mentioned before, there is still much to be discovered about the multivariate Gaussian distribution. In particular, we found comparisons of moments to be scarce in the literature even for the bivariate case. In an attempt to address this gap, we discovered a very close relationship between Gaussian moment ratio inequalities (MRIs) and the GPI. On one hand, in \cite[Theorem 4.3]{RusSun} we gave a novel bivariate MRI that was weaker than the GPI. On the other hand, in this paper, we prove a new MRI (see Theorem \ref{thm2} below) that is stronger than Theorem \ref{thmm1m2m3}.

Define
$$
{\cal S}:=\{(1,m_3):m_3\ge 5\}\bigcup\{(2,m_3):m_3\ge 3\}\bigcup\{(m_2,m_3):m_3\ge m_2\ge3\},
$$
$$
r_{m_2,m_3}=(2m_2+1)(2m_3+1)+1,\ \ \ \ t_{m_2,m_3}=\frac{1}{r_{m_2,m_3}+\left(1+\frac{1}{2m_2}\right)\left(1+\frac{1}{2m_3}\right)},
$$
and for $\frac{1}{r^2_{m_2,m_3}}<z\le 1$,
\begin{eqnarray*}
&&H_{m_2,m_3}(z)\\
&=&\frac{(m_2+m_3+1)(r_{m_2,m_3}z-1)+\sqrt{[(m_3-m_2)(r_{m_2,m_3}z-1)]^2+(2m_2+1)^3(2m_3+1)^3z}}{r_{m_2,m_3}^2z-1}.
\end{eqnarray*}
For random variables $X$ and $Y$,  denote by ${\rm Cov}(X,Y)$ and ${\rm Corr}(X,Y)$ their covariance and correlation, respectively. We will show that the GPI (\ref{GPIk}) is implied by the following MRI.

\begin{thm}\label{thm2}
Let $(X_2,X_3)$ be a centered Gaussian random vector. If $(m_2,m_3)\in{\cal S}$, then
\begin{eqnarray}\label{MRI22}
&&\frac{\left|E[ X_2^{2m_2+1}X_3^{2m_3+1}]\right|}{(2m_2+1)(2m_3+1)E[ X_2^{2m_2}X_3^{2m_3}]}\nonumber\\
&\le&
 \begin{cases}
 |{\rm Cov}(X_2,X_3)|,\ \  &{\rm if}\ |{\rm Corr}(X_2,X_3)|\le \sqrt{t_{m_2,m_3}},\nonumber\\
H_{m_2,m_3}([{\rm Corr}(X_2,X_3)]^2)\cdot|{\rm Cov}(X_2,X_3)|,\ \ &{\rm if}\ \sqrt{t_{m_2,m_3}}<|{\rm Corr}(X_2,X_3)|.
   \end{cases}\nonumber\\
&&
\end{eqnarray}
The equality sign holds if and only $X_2$ and $X_3$ are independent.
\end{thm}

\noindent Note that the MRI (\ref{MRI22}) does not hold if $m_2=1$ and $1\le m_3\le 4$, or $m_2=m_3=2$. However, for these particular cases, the  GPI (\ref{GPIk}) has been proved by \cite{LHS}. These cases can also be easily handled  by the SOS method developed in \cite{RusSun2}.

In addition to the intrinsic value of solving the GPI and the implications a proof will have on related problems and fields, we can see that the study of the GPI itself has already led to some interesting new results, including the MRIs discussed above, and highly non-trivial combinatorial identities such as \cite[Lemma 2.5 and Corollary 2.8]{LHS}. Therefore, the GPI certainly deserves careful consideration. Finally, we hope that the interplay between computing and hard analysis that we  describe in this paper will inspire confidence in the research community to extend our results.

The remainder of this paper is structured as follows. In Section \ref{sec2}, we explain the relation between the GPI (\ref{GPIk}) and the MRI (\ref{MRI22}). First, we derive an  inequality (see (\ref{Augsdyy}) below) which is slightly stronger than the GPI. Then, we show that this  inequality is equivalent to a hypergeometric function ratio inequality (HFRI)  (see (\ref{Aug18a}) below) and hence is equivalent to the MRI (\ref{MRI22}). In Sections \ref{sec3} and \ref{sec4},  we use different methods to prove the HFRI (\ref{Aug18a}) for cases $1\le m_2\le 7$ and $m_2\ge 8$ separately.  In Section \ref{sec5}, we make some concluding remarks. To complete the proof of Section \ref{sec4}, we need the sums-of-squares (SOS) expansion of a difference function. This expansion is obtained by \textbf{Mathematica} and given in Section \ref{sec6}, the Appendix.

To simplify notation, in the sequel, we use $r$, $t$ and $H$ to denote $r_{m_2,m_3}$, $t_{m_2,m_3}$ and $H_{m_2,m_3}$, respectively, whenever no confusion is caused. Throughout this paper, we assume without loss of generality that $m_3\ge m_2$.

\section{Relation between GPI and MRI}\label{sec2}\setcounter{equation}{0}

In this section, we derive the MRI (\ref{MRI22}) and show that its validity implies  the validity of the GPI (\ref{GPIk}).

Denote by $F(a, b; c; z)$ the hypergeometric function (cf. \cite{R}):
$$
    F(a, b; c; z)=\sum_{j=0}^{\infty}\frac{(a)_{j}(b)_{j}}{(c)_{j}}\cdot\frac{z^j}{j!},\ \ \ \ |z|<1,
$$
where $(\alpha)_{j}:=\alpha(\alpha+1)\cdots(\alpha+j-1)$ for $j\ge 1$, and $(\alpha)_0=1$ for $\alpha\not=0$. Let $(X_2,X_3)$ be a centered Gaussian random vector and $x={\rm Corr}(X_2,X_3)$. By the moment formula (cf. \cite[Page 261]{KBJ}), we get
\begin{eqnarray}\label{we2}
&&E[ X_2^{2m_2}X_3^{2m_3}]\nonumber\\
&=&(2m_2-1)!!(2m_3-1)!![{\rm Var}(X_2)]^{m_2}[{\rm Var}(X_3)]^{m_3}F\left(-m_2,-m_3;\frac{1}{2};x^2\right),
\end{eqnarray}
and
\begin{eqnarray}\label{we1}
&&E[ X_2^{2m_2+1}X_3^{2m_3+1}]\nonumber\\
&=&(2m_2+1)!!(2m_3+1)!![{\rm Var}(X_2)]^{\frac{2m_2+1}{2}}[{\rm Var}(X_3)]^{\frac{2m_3+1}{2}}xF\left(-m_2,-m_3;\frac{3}{2};x^2\right).\ \ \ \
\end{eqnarray}

Consider the following functions contiguous to ${F}(a,b;c;z)$:
\begin{equation}\label{contiguous}
  {F}(a\pm 1,b;c;z),\quad {F}(a,b\pm 1; c;z),\quad {F}(a,b;c\pm 1;z).
\end{equation}
To simplify notation,  we denote ${F}(a,b;c;z)$ and
the six  functions in (\ref{contiguous}) respectively by
$$
  {F},\quad {F}(a\pm 1),\quad {F}(b\pm 1),\quad {F}(c\pm 1).
$$
We have the following relations of Gauss between contiguous functions (cf. \cite[2.8-(21), (31), (37), (38), pages 102 and 103]{B}):
\begin{eqnarray}\label{Gauss3}
&&F'=\frac{a[F(a+1)-F]}{z},\nonumber\\
&&[c-2a-(b-a)z]F+a(1-z)F(a+1)-(c-a)F(a-1)=0,\nonumber\\
&&(b-a)(1-z)F-(c-a)F(a-1)+(c-b)F(b-1)=0,\nonumber\\
&&c(1-z)F-cF(a-1)+(c-b)zF(c+1)=0.
\end{eqnarray}

\subsection{A stronger inequality}

By \cite[Lemma 2.1]{RusSun}, to prove the GPI (\ref{GPIk}), we may assume without loss of generality that $X_1 = X_2+aX_3$ for some $a \in \mathbb{R}$ and $E[X_2^2]=E[X_3^2]=1$.
Define
$$x=E[X_2 X_3].
$$
Then,
$$
E[X_1^2]=a^2+1+2ax.
$$

Consider the moment ratio:
\begin{eqnarray}\label{moment}
&&\frac{E[(X_2+aX_3)^{2}X_2^{2m_2} X_3^{2m_3} ]} {E[X_1^{2}] E[X_2^{2m_2}] E[X_3^{m_3}]}\nonumber\\
&=&\frac{E[a^2X_2^{2m_2} X_3^{2m_3+2} + X_2^{2m_2+2} X_3^{2m_3} + 2aX_2^{2m_2+1} X_3^{2m_2+1}]} {(a^2+1+2ax)(2m_2-1)!!(2m_3-1)!!}\nonumber\\
&=&\Bigg[a^2\left(2m_3+1\right)F\left(-m_3-1,-m_2;\frac{1}{2};x^2\right)+\left(2m_2+1\right)F\left(-m_3,-m_2-1;\frac{1}{2};x^2\right)\nonumber\\
&&+2ax\left(2m_3+1\right)\left(2m_2+1\right) F\left(-m_3,-m_2;\frac{3}{2};x^2\right)\Bigg]\cdot\frac{1}{a^2+1+2ax}.
\end{eqnarray}
Then, the proof of Theorem \ref{thmm1m2m3} is complete if we can show that for any $a\in\mathbb{R}$ and $x\in[-1,1]$,
\begin{eqnarray}\label{gen2}
&&a^2\left(2m_3+1\right)F\left(-m_3-1,-m_2;\frac{1}{2};x^2\right)+\left(2m_2+1\right)F\left(-m_3,-m_2-1;\frac{1}{2};x^2\right)\nonumber\\
&&+2ax\left(2m_3+1\right)\left(2m_2+1\right) F\left(-m_3,-m_2;\frac{3}{2};x^2\right)\nonumber\\
&>&a^2+1+2ax.
\end{eqnarray}
Equivalently, we need to show that for any $a\in(-\infty,0]$ and $x\in[0,1]$,
\begin{eqnarray}\label{Aug9}
F(a,x)&:=&a^2\left[\left(2m_3+1\right)F\left(-m_2,-m_3-1;\frac{1}{2};x^2\right)-1\right]\nonumber\\
&&+2ax\left[\left(2m_3+1\right)\left(2m_2+1\right) F\left(-m_2,-m_3;\frac{3}{2};x^2\right)-1\right]\nonumber\\
&&+\left[\left(2m_2+1\right)F\left(-m_2-1,-m_3;\frac{1}{2};x^2\right)-1\right]\nonumber\\
&>&0.
\end{eqnarray}
Note that (\ref{Aug9}) holds when either $x\in\{0,1\}$ or $a=0$, and
$$
\lim_{a\rightarrow-\infty}\min_{x\in[0,1]}F(a,x)=\infty.
$$
Thus, we need only show that (\ref{Aug9}) holds under the following condition:
$$
\frac{\partial F(a,x)}{\partial a}=0,\ \ \ \ a\in(-\infty,0),\, x\in(0,1).
$$

By (\ref{Gauss3}), we get
\begin{eqnarray*}
&&F\left(-m_2,-m_3-1;\frac{1}{2};x^2\right)\\
&=&\frac{1}{2m_3+1}\left[(2m_2+1)F\left(-m_2-1,-m_3;\frac{1}{2};x^2\right)+2(m_3-m_2)(1-x^2)F\left(-m_2,-m_3;\frac{1}{2};x^2\right)\right],
\end{eqnarray*}
and
\begin{eqnarray}\label{Aug2222}
&&x^2F\left(-m_2,-m_3;\frac{3}{2};x^2\right)\nonumber\\
&=&\frac{1}{2m_3+1}\left[F\left(-m_2-1,-m_3;\frac{1}{2};x^2\right)-(1-x^2)F\left(-m_2,-m_3;\frac{1}{2};x^2\right)\right].
\end{eqnarray}
Then,
\begin{eqnarray*}
&&\frac{\partial F(a,x)}{\partial a}=0\nonumber\\
&\Leftrightarrow&0=a\left[\left(2m_3+1\right)F\left(-m_2,-m_3-1;\frac{1}{2};x^2\right)-1\right]\nonumber\\
&&\ \ \ \ \ +x\left[\left(2m_3+1\right)\left(2m_2+1\right) F\left(-m_2,-m_3;\frac{3}{2};x^2\right)-1\right]\nonumber\\
&\Leftrightarrow&0=ax\left[(2m_2+1)F\left(-m_2-1,-m_3;\frac{1}{2};x^2\right)+2(m_3-m_2)(1-x^2)F\left(-m_2,-m_3;\frac{1}{2};x^2\right)-1\right]\nonumber\\
&&\ \ \ \ \ +\left\{\left(2m_2+1\right)\left[F\left(-m_2-1,-m_3;\frac{1}{2};x^2\right)-(1-x^2)F\left(-m_2,-m_3;\frac{1}{2};x^2\right)\right]-x^2\right\}\nonumber\\
&\Leftrightarrow&a=-\frac{\left(2m_2+1\right)\left[F\left(-m_2-1,-m_3;\frac{1}{2};x^2\right)-(1-x^2)F\left(-m_2,-m_3;\frac{1}{2};x^2\right)\right]-x^2}{x\left[(2m_2+1)F\left(-m_2-1,-m_3;\frac{1}{2};x^2\right)+2(m_3-m_2)(1-x^2)F\left(-m_2,-m_3;\frac{1}{2};x^2\right)-1\right]}.\nonumber\\
&&
\end{eqnarray*}
Hence, we have that
\begin{eqnarray*}
F(a,x)&=&a^2\left[(2m_2+1)F\left(-m_2-1,-m_3;\frac{1}{2};x^2\right)+2(m_3-m_2)(1-x^2)F\left(-m_2,-m_3;\frac{1}{2};x^2\right)-1\right]\nonumber\\
&&+\frac{2a}{x}\left\{\left(2m_2+1\right) \left[F\left(-m_2-1,-m_3;\frac{1}{2};x^2\right)-(1-x^2)F\left(-m_2,-m_3;\frac{1}{2};x^2\right)\right]-x^2\right\}\nonumber\\
&&+\left[\left(2m_2+1\right)F\left(-m_2-1,-m_3;\frac{1}{2};x^2\right)-1\right]\nonumber\\
&=&\frac{\{\left(2m_2+1\right)\left[F\left(-m_2-1,-m_3;\frac{1}{2};x^2\right)-(1-x^2)F\left(-m_2,-m_3;\frac{1}{2};x^2\right)\right]-x^2\}^2}{x^2\left[(2m_2+1)F\left(-m_2-1,-m_3;\frac{1}{2};x^2\right)+2(m_3-m_2)(1-x^2)F\left(-m_2,-m_3;\frac{1}{2};x^2\right)-1\right]}\nonumber\\
&&-\frac{2\{\left(2m_2+1\right)\left[F\left(-m_2-1,-m_3;\frac{1}{2};x^2\right)-(1-x^2)F\left(-m_2,-m_3;\frac{1}{2};x^2\right)\right]-x^2\}^2}{x^2\left[(2m_2+1)F\left(-m_2-1,-m_3;\frac{1}{2};x^2\right)+2(m_3-m_2)(1-x^2)F\left(-m_2,-m_3;\frac{1}{2};x^2\right)-1\right]}\nonumber\\
&&+\left[\left(2m_2+1\right)F\left(-m_2-1,-m_3;\frac{1}{2};x^2\right)-1\right].
\end{eqnarray*}
Thus, to prove the GPI (\ref{GPIk}), we need to show that
\begin{eqnarray}\label{Aug111}
&&\left\{\left(2m_2+1\right)\left[F\left(-m_2-1,-m_3;\frac{1}{2};x^2\right)-(1-x^2)F\left(-m_2,-m_3;\frac{1}{2};x^2\right)\right]-x^2\right\}^2\nonumber\\
&<&x^2\left[\left(2m_2+1\right)F\left(-m_2-1,-m_3;\frac{1}{2};x^2\right)-1\right]\nonumber\\
&&\cdot\left[(2m_2+1)F\left(-m_2-1,-m_3;\frac{1}{2};x^2\right)+2(m_3-m_2)(1-x^2)F\left(-m_2,-m_3;\frac{1}{2};x^2\right)-1\right].\nonumber\\
&&
\end{eqnarray}
It can be shown that (\ref{Aug111}) is equivalent to \cite[Claim C, page 14]{RusSun}.

By (\ref{Gauss3}), we get
\begin{eqnarray*}
&&(1-x^2)F\left(-m_2,-m_3;\frac{1}{2};x^2\right)\\
&=&F\left(-m_2-1,-m_3;\frac{1}{2};x^2\right)-(2m_3+1)x^2F\left(-m_2,-m_3;\frac{3}{2};x^2\right),
\end{eqnarray*}
and
\begin{eqnarray*}
0&=&-(m_3-m_2)(1-x^2)F\left(-m_2,-m_3;\frac{1}{2};x^2\right)\\
&&-\left(m_2+\frac{1}{2}\right)F\left(-m_2-1,-m_3;\frac{1}{2};x^2\right)+\left(m_3+\frac{1}{2}\right)F\left(-m_2,-m_3-1;\frac{1}{2};x^2\right).
\end{eqnarray*}
Then, (\ref{Aug111}) becomes
\begin{eqnarray}\label{we3}
&&x^2\left\{\left(2m_2+1\right)(2m_3+1)F\left(-m_2,-m_3;\frac{3}{2};x^2\right)-1\right\}^2\nonumber\\
&<&\left[\left(2m_2+1\right)F\left(-m_2-1,-m_3;\frac{1}{2};x^2\right)-1\right]\nonumber\\
&&\cdot\left[\left(2m_3+1\right)F\left(-m_2,-m_3-1;\frac{1}{2};x^2\right)-1\right].
\end{eqnarray}
Note that, different from (\ref{Aug111}), the roles played by $m_2$ and $m_3$ in (\ref{we3}) are symmetric. This key fact will lead us to derive the stronger inequality (\ref{Augsdyy}) given below.

\begin{rem}
Let $(X_2,X_3)$ be a centered Gaussian random vector with ${\rm Var}(X_2)={\rm Var}(X_3)=1$. By (\ref{we2}) and (\ref{we1}), we can rewrite (\ref{we3}) as the following interesting and delicate moment inequality, which is equivalent to the GPI (\ref{GPIk}).
\begin{eqnarray*}
&&\left\{E[ X_2^{2m_2+1}X_3^{2m_3+1}]-E[X_2^{2m_2}]E[X_3^{2m_3}]E[X_2X_3]\right\}^2\nonumber\\
&<&\left\{E[X_2^{2m_2+2}X_3^{2m_3}]-E[X_2^{2m_2}]E[X_3^{2m_3}]\right\}\\
&&\cdot\left\{E[X_2^{2m_2}X_3^{2m_3+2}]-E[X_2^{2m_2}]E[X_3^{2m_3}]\right\}.
\end{eqnarray*}
\end{rem}

Let $z=x^2$. Then, (\ref{we3}) becomes
\begin{eqnarray}\label{Aug333}
&&z\left\{\left(2m_2+1\right)(2m_3+1)F\left(-m_2,-m_3;\frac{3}{2};z\right)-1\right\}^2\nonumber\\
&<&\left[\left(2m_2+1\right)F\left(-m_2-1,-m_3;\frac{1}{2};z\right)-1\right]\nonumber\\
&&\cdot\left[\left(2m_3+1\right)F\left(-m_2,-m_3-1;\frac{1}{2};z\right)-1\right].
\end{eqnarray}
Further, by (\ref{Aug2222}), we find that (\ref{Aug333}) is equivalent to
\begin{eqnarray*}
&&z\left\{\left(2m_2+1\right)(2m_3+1)F\left(-m_2,-m_3;\frac{3}{2};z\right)-1\right\}^2\\
&<&\left[\left(2m_2+1\right)(1-z)F\left(-m_2,-m_3;\frac{1}{2};z\right)+\left(2m_2+1\right)(2m_3+1)zF\left(-m_2,-m_3;\frac{3}{2};z\right)-1\right]\\
&&\cdot\left[\left(2m_3+1\right)(1-z)F\left(-m_2,-m_3;\frac{1}{2};z\right)+\left(2m_2+1\right)(2m_3+1)zF\left(-m_2,-m_3;\frac{3}{2};z\right)-1\right],
\end{eqnarray*}
i.e.,
\begin{eqnarray*}
&&z\left\{\left(2m_2+1\right)(2m_3+1)F\left(-m_2,-m_3;\frac{3}{2};z\right)-1\right\}^2\\
&<&\left[\left\{\left(2m_2+1\right)F\left(-m_2,-m_3;\frac{1}{2};z\right)-1\right\}(1-z)+\left(2m_2+1\right)(2m_3+1)zF\left(-m_2,-m_3;\frac{3}{2};z\right)-z\right]\\
&&\cdot\left[\left\{\left(2m_3+1\right)F\left(-m_2,-m_3;\frac{1}{2};z\right)-1\right\}(1-z)+\left(2m_2+1\right)(2m_3+1)zF\left(-m_2,-m_3;\frac{3}{2};z\right)-z\right].
\end{eqnarray*}
Hence, we need only show that for any $z\in(0,1)$,
\begin{eqnarray}\label{Augsd}
z
&<&\left[\frac{\left\{\left(2m_2+1\right)F\left(-m_2,-m_3;\frac{1}{2};z\right)-1\right\}(1-z)}{\left(2m_2+1\right)(2m_3+1)F\left(-m_2,-m_3;\frac{3}{2};z\right)-1}+z\right]\nonumber\\
&&\cdot\left[\frac{\left\{\left(2m_3+1\right)F\left(-m_2,-m_3;\frac{1}{2};z\right)-1\right\}(1-z)}{\left(2m_2+1\right)(2m_3+1)F\left(-m_2,-m_3;\frac{3}{2};z\right)-1}+z\right].
\end{eqnarray}

Note that
\begin{eqnarray*}
&&\frac{\left\{\left(2m_2+1\right)F\left(-m_2,-m_3;\frac{1}{2};z\right)-1\right\}(1-z)}{\left(2m_2+1\right)(2m_3+1)F\left(-m_2,-m_3;\frac{3}{2};z\right)-1}\\
&\ge&\frac{\left(2m_2+1\right)F\left(-m_2,-m_3;\frac{1}{2};z\right)(1-z)}{\left(2m_2+1\right)(2m_3+1)F\left(-m_2,-m_3;\frac{3}{2};z\right)}+\frac{-(1-z)}{\left(2m_2+1\right)(2m_3+1)F\left(-m_2,-m_3;\frac{3}{2};z\right)}\\
&\ge&\frac{F\left(-m_2,-m_3;\frac{1}{2};z\right)(1-z)}{(2m_3+1)F\left(-m_2,-m_3;\frac{3}{2};z\right)}+\frac{-(1-z)}{\left(2m_2+1\right)(2m_3+1)}.
\end{eqnarray*}
Therefore, to prove (\ref{Augsd}), it suffices to show that the following stronger inequality holds for  $0<z<1$:
\begin{eqnarray}\label{Augsdyy}
z
&<&\left[\frac{F\left(-m_2,-m_3;\frac{1}{2};z\right)(1-z)}{(2m_3+1)F\left(-m_2,-m_3;\frac{3}{2};z\right)}+\frac{[(2m_2+1)(2m_3+1)+1]z-1}{\left(2m_2+1\right)(2m_3+1)}\right]\nonumber\\
&&\cdot\left[\frac{F\left(-m_2,-m_3;\frac{1}{2};z\right)(1-z)}{(2m_2+1)F\left(-m_2,-m_3;\frac{3}{2};z\right)}+\frac{[(2m_2+1)(2m_3+1)+1]z-1}{\left(2m_2+1\right)(2m_3+1)}\right].
\end{eqnarray}
Using \textbf{Mathematica}, it appears that inequality (\ref{Augsdyy}) holds if $
(m_2,m_3)\in{\cal S}$. In Sections \ref{sec3} and \ref{sec4} below, we will give a rigorous proof for  (\ref{Augsdyy}) by combining  computing and hard analysis.

\subsection{MRI $\Leftrightarrow$ HFRI $\Rightarrow$ GPI}

For $0<z< 1$, define
$$
y=\frac{F\left(-m_2,-m_3;\frac{1}{2};z\right)}{F\left(-m_2,-m_3;\frac{3}{2};z\right)}.
$$
Then, $y>1$ and
\begin{eqnarray*}
&&(\ref{Augsdyy})\ {\rm holds}\nonumber\\
&\Leftrightarrow&0<\frac{(1-z)^2}{(2m_2+1)(2m_3+1)}\cdot y^2\nonumber\\
&&\ \ \ \ +\frac{2(m_2+m_3+1)\{[(2m_2+1)(2m_3+1)+1]z-1\}(1-z)}{(2m_2+1)^2(2m_3+1)^2}\cdot y\nonumber\\
&&\ \ \ \ +\frac{[(2m_2+1)(2m_3+1)+1]^2z^2-\{[(2m_2+1)(2m_3+1)+1]^2+1\}z+1}{(2m_2+1)^2(2m_3+1)^2}\nonumber\\
&\Leftrightarrow&0<(1-z)y^2+2\beta y+\gamma,
\end{eqnarray*}
where
\begin{eqnarray*}
&&\beta:=-\frac{(m_2+m_3+1)\{1-[(2m_2+1)(2m_3+1)+1]z\}}{(2m_2+1)(2m_3+1)},\nonumber\\
&&\gamma:=\frac{1-[(2m_2+1)(2m_3+1)+1]^2z}{(2m_2+1)(2m_3+1)}.
\end{eqnarray*}

By the identity
$$
(m_2+m_3+1)^2=(m_3-m_2)^2+(2m_2+1)(2m_3+1),
$$
we get
\begin{eqnarray}\label{KKJJGG}
&&\beta^2-(1-z)\gamma\nonumber\\
&=&\left(\frac{(m_3-m_2)\{1-[(2m_2+1)(2m_3+1)+1]z\}}{(2m_2+1)(2m_3+1)}\right)^2+(2m_2+1)(2m_3+1)z.
\end{eqnarray}
Then,
\begin{eqnarray*}
&&0<(1-z)y^2+2\beta y+\gamma\nonumber\\
&\Leftrightarrow&\left|y+\frac{\beta}{1-z}\right|>\frac{\sqrt{\beta^2-(1-z)\gamma}}{1-z}.
\end{eqnarray*}
For $0<z\le\frac{1}{r^2}$, we have
\begin{eqnarray*}
&&y+\frac{\beta}{1-z}-\frac{\sqrt{\beta^2-(1-z)\gamma}}{1-z}\nonumber\\
&>&1-\frac{m_2+m_3+1}{(2m_2+1)(2m_3+1)\left[1-\frac{1}{[(2m_2+1)(2m_3+1)+1]^2}\right]}\nonumber\\
&&-\frac{\sqrt{\left[\frac{m_3-m_2}{(2m_2+1)(2m_3+1)}\right]^2+\frac{(2m_2+1)(2m_3+1)}{[(2m_2+1)(2m_3+1)+1]^2}}}{1-\frac{1}{[(2m_2+1)(2m_3+1)+1]^2}}\nonumber\\
&>&1-\frac{1}{(2m_2+1)\left(1-\frac{1}{10^2}\right)}-\frac{\sqrt{\left[\frac{1}{2(2m_2+1)}\right]^2+\frac{1}{(2m_2+1)^2}}}{1-\frac{1}{10^2}}\nonumber\\
&>&0.
\end{eqnarray*}
Hence, to prove (\ref{Augsdyy}), it suffices to show that for $\frac{1}{r^2}< z<1$,
$$
y>\frac{-\beta+\sqrt{\beta^2-(1-z)\gamma}}{1-z}.
$$

Note that
\begin{eqnarray}\label{HGF}
\frac{1}{H(z)}&=&\frac{\gamma}{-\beta-\sqrt{\beta^2-(1-z)\gamma}}\nonumber\\
&=&\frac{-\beta+\sqrt{\beta^2-(1-z)\gamma}}{1-z},\ \ \ \ 0<z<1.
\end{eqnarray}
Hence we need show that the following HFRI holds:
\begin{eqnarray*}
\frac{F\left(-m_2,-m_3;\frac{1}{2};z\right)}{F\left(-m_2,-m_3;\frac{3}{2};z\right)}>\frac{1}{H(z)},\ \ \ \ \frac{1}{r^2}< z<1.
\end{eqnarray*}
For  $\frac{1}{r^2}< z<\frac{1}{r}$, we have
\begin{eqnarray*}
&&H(z)\ge1\nonumber\\
&\Leftrightarrow&[(m_3-m_2)(rz-1)]^2+(r-1)^3z\ge[(r^2z-1)+(m_2+m_3+1)(1-rz)]^2\nonumber\\
&\Leftrightarrow&(r-1)^3z-(r^2z-1)^2-(r-1)(rz-1)^2\ge2(m_2+m_3+1)(1-rz)(r^2z-1)\nonumber\\
&\Leftrightarrow&(r-1)(1-z)(r^2z-1)-(r^2z-1)^2\ge2(m_2+m_3+1)(1-rz)(r^2z-1)\nonumber\\
&\Leftrightarrow&(r-1)(1-z)-(r^2z-1)\ge2(m_2+m_3+1)(1-rz)\nonumber\\
&\Leftrightarrow&r(1-rz)-(r-1)z\ge2(m_2+m_3+1)(1-rz)\nonumber\\
&\Leftrightarrow&(2m_2)(2m_3)(1-rz)\ge(r-1)z\nonumber\\
&\Leftrightarrow&z\le t.
\end{eqnarray*}
Then,
\begin{eqnarray*}
\frac{F\left(-m_2,-m_3;\frac{1}{2};z\right)}{F\left(-m_2,-m_3;\frac{3}{2};z\right)}>\frac{1}{H(z)},\ \ \ \ \frac{1}{r^2}<z\le t.
\end{eqnarray*}
Thus, to prove the GPI (\ref{GPIk}), we need only prove
\begin{eqnarray}\label{Aug18a}
\frac{F\left(-m_2,-m_3;\frac{1}{2};z\right)}{F\left(-m_2,-m_3;\frac{3}{2};z\right)}>\frac{1}{H(z)},\ \ \ \ t<z< 1.
\end{eqnarray}

Finally, we show that the HFRI (\ref{Aug18a}) is equivalent to the MRI (\ref{MRI22}) with the equality sign holding if and only if $X_2$ and $X_3$ are independent. We assume without loss of generality that ${\rm Var}(X_2)={\rm Var}(X_3)=1$.  By (\ref{we2}) and (\ref{we1}), we get
\begin{eqnarray}\label{Aug27a}
\frac{|E[ X_2^{2m_2+1}X_3^{2m_3+1}]|}{E[ X_2^{2m_2}X_3^{2m_3}]}={(2m_2+1)(2m_3+1)|x|}\frac{F\left(-m_2,-m_3;\frac{3}{2};x^2\right)}{F\left(-m_2,-m_3;\frac{1}{2};x^2\right)}.
\end{eqnarray}
Note that
$$
\frac{F\left(-m_2,-m_3;\frac{1}{2};0\right)}{F\left(-m_2,-m_3;\frac{3}{2};0\right)}=1,
$$
and
$$
\frac{F\left(-m_2,-m_3;\frac{1}{2};z\right)}{F\left(-m_2,-m_3;\frac{3}{2};z\right)}>1,\ \ \ \ 0<z< 1.
$$
Therefore, by  (\ref{Aug27a}) and (\ref{HGFD}) (see below), we conclude that  (\ref{Aug18a}) is equivalent to (\ref{MRI22}) with the equality sign holding if and only if $X_2$ and $X_3$ are independent.

\begin{rem}
On one hand, the MRI (\ref{MRI22}) implies the GPI (\ref{GPIk}). On the other hand, \cite[Theorem 4.3 and Remark 4.4]{RusSun} shows that  the following MRI holds and is implied by the GPI (\ref{GPIk}):
\begin{eqnarray}\label{MRI33}
\frac{\left|E[ X_2^{2m_2-1}X_3^{2m_3+1}]\right|}{E[ X_2^{2m_2}X_3^{2m_3}]}<\frac{2m_3+1}{2m_2},
\end{eqnarray}
where $(X_2,X_3)$ is a centered Gaussian random vector with ${\rm Var}(X_2)={\rm Var}(X_3)=1$. Hence the MRI (\ref{MRI22}) is stronger than the MRI (\ref{MRI33}). However, these two MRIs have independent interests.
\end{rem}

\section{Proof of HFRI (\ref{Aug18a}) for the case $1\le m_2\le 7$}\label{sec3}\setcounter{equation}{0}

In this section, we show that (\ref{Aug18a}) holds for $1\le m_2\le 7$.

\subsection{Positiveness of bivariate polynomials}

Note that
\begin{eqnarray}\label{Aug27HH}
F\left(-m_2,-m_3;\frac{1}{2};z\right)&=&m_2!m_3!\sum_{j=0}^{m_2}\frac{2^{2j}z^j}{(m_2-j)!(m_3-j)!(2j)!},\nonumber\\
F\left(-m_2,-m_3;\frac{3}{2};z\right)&=&m_2!m_3!\sum_{j=0}^{m_2}\frac{2^{2j}z^j}{(m_2-j)!(m_3-j)!(2j+1)!}.
\end{eqnarray}
To simplify notation, we denote
$$
f_1(z):=F\left(-m_2,-m_3;\frac{1}{2};z\right),\ \ \ \ f_2(z):=F\left(-m_2,-m_3;\frac{3}{2};z\right).
$$
Then, we have
\begin{eqnarray}\label{Aug22g}
&&(\ref{Aug18a})\ {\rm holds}\nonumber\\
&\Leftarrow&\ \ \left[(m_2+m_3+1)(rz-1)+\sqrt{[(m_3-m_2)(rz-1)]^2+(2m_2+1)^3(2m_3+1)^3z}\right]f_1(z)\nonumber\\
&&>(r^2z-1)f_2(z)\nonumber\\
&\Leftarrow&\ \ \left\{[(m_3-m_2)(rz-1)]^2+(r-1)^3z\right\}f^2_1(z)\nonumber\\
&&>\left\{(r^2z-1)f_2(z)-(m_2+m_3+1)(rz-1)f_1(z)\right\}^2\nonumber\\
&\Leftrightarrow&\ \ (r-1)^3zf^2_1(z)\nonumber\\
&&>(r^2z-1)^2f_2^2(z)+(r-1)(rz-1)^2f_1^2(z)-2(m_2+m_3+1)(rz-1)(r^2z-1)f_1(z)f_2(z)\nonumber\\
&\Leftrightarrow&\ \ (r-1)(1-z)(r^2z-1)f_1^2(z)\nonumber\\
&&>(r^2z-1)^2f_2^2(z)-2(m_2+m_3+1)(rz-1)(r^2z-1)f_1(z)f_2(z)\nonumber\\
&\Leftrightarrow&\ \ (r-1)(1-z)f_1^2(z)+2(m_2+m_3+1)(rz-1)f_1(z)f_2(z)\nonumber\\
&&>(r^2z-1)f_2^2(z)\nonumber\\
&\Leftrightarrow&S_{m_2, m_3}(z):=(2m_2+1)(2m_3+1)\left(1-z\right)\left[\sum_{j=0}^{m_2}\frac{2^{2j}z^jm_2!m_3!}{(m_2-j)!(m_3-j)!(2j)!}\right]^2\nonumber\\
&&\ \ \ \ \ \ \ \ \ \ \ \ \ \ \ \ \ +2(m_2+m_3+1)\left\{{[(2m_2+1)(2m_3+1)+1]z}-1\right\}\nonumber\\
&&\ \ \ \ \ \ \ \ \ \ \ \ \ \ \ \ \ \ \ \cdot\left[\sum_{j=0}^{m_2}\frac{2^{2j}z^jm_2!m_3!}{(m_2-j)!(m_3-j)!(2j)!}\right]\left[\sum_{j=0}^{m_2}\frac{2^{2j}z^jm_2!m_3!}{(m_2-j)!(m_3-j)!(2j+1)!}\right]\nonumber\\
&&\ \ \ \ \ \ \ \ \ \ \ \ \ \ \ \ \ -\left\{{[(2m_2+1)(2m_3+1)+1]^2z}-1\right\}\left[\sum_{j=0}^{m_2}\frac{2^{2j}z^jm_2!m_3!}{(m_2-j)!(m_3-j)!(2j+1)!}\right]^2\nonumber\\
&&\ >0,\ \ \ \ \frac{1}{r^2}<z<1.
\end{eqnarray}

By using the transformations
$$
z=\frac{c^2}{1+c^2},
$$
and
\begin{eqnarray*}
m_3=
 \begin{cases}
b^2+5,\ \  &{\rm if}\ m_2=1,\\
b^2+3,\ \  &{\rm if}\ m_2=2,\\
b^2+m_2,\ \  &{\rm if}\ m_2\ge3,
   \end{cases}
\end{eqnarray*}
we define
$$
h_{m_2}(b,c):=(1+c^2)^{2m_2+1}\cdot S_{m_2, m_3}(z),\ \ \ \ b,c \in \mathbb{R}.
$$
Then, by (\ref{Aug22g}), we find that the HFRI (\ref{Aug18a}) holds for $1\le m_2\le 7$ if  $h_{m_2}$ is a positive bivariate polynomial for each $1\le m_2\le 7$.

By virtue of \textbf{Mathematica}, we obtain the expansion of $h_{m_2}$, a polynomial of $b$ and $c$. See the next subsection. Note that all exponents of $b$ and $c$ are even. We sum all terms with negative coefficients with some of the terms with positive coefficients and then apply the package `\texttt{SumsOfSquares}' in \textbf{Macaulay2} \cite{SOS,PP} to obtain exact SOS decompositions. Hence, since $h_{m_2}$ is strictly positive, we conclude that
\begin{eqnarray*}
(\ref{Aug18a})\ {\rm holds\ for}\ 1\le m_2\le 7.
\end{eqnarray*}

\subsection{Expansion of functions $h_{m_2}$}

Let $1\le m_2\le 7$. We have the following expansions of the functions $h_{m_2}$ for any $b,c \in \mathbb{R}$ and the SOS decompositions corresponding to the sums of the 10 terms in parentheses:

\begin{dmath*}
h_1(b,c)=\frac{8 b^6 c^6}{3}+\frac{124 b^4 c^6}{3}+\frac{622 b^2
   c^6}{3}+\frac{1}{9} \left(48 b^6 c^4+664 b^4 c^4-48 b^4 c^2+2884 b^2
   c^4-546 b^2 c^2+36 b^2+3003 c^6+3766 c^4-1557 c^2+180\right).
\end{dmath*}
Using `\texttt{SumsOfSquares}' in \textbf{Macaulay2}, we get
\begin{dmath*}
48 b^6 c^4+664 b^4 c^4-48 b^4 c^2+2884 b^2 c^4-546 b^2 c^2+36 b^2+3003 c^6+3766 c^4-1557 c^2+180 = 4924\bigg(\frac{3377}{39392}\,b^{2}c^{2}+c^{2}-\frac{35173}{196960} \bigg)^2
+ 3003\bigg(\frac{151}{8008}\,b^{2}c+c^{3}-\frac{193}{1001}\,c \bigg)^2
+\frac{3853}{2}\bigg(\frac{557}{77060}\,b^{3}c^{2}+b\,c^{2}-\frac{16027}{200356}\,b \bigg)^2
+\frac{945348187}{1575680}\bigg(b^{2}c^{2}-\frac{1722648687}{12289526431} \bigg)^2
+\frac{1591}{13}\bigg(b\,c \bigg)^2
+\frac{1802093}{20020}\bigg(-\frac{193955}{1802093}\,b^{2}c+c \bigg)^2
+\frac{147644951}{3082400}\bigg(b^{3}c^{2}-\frac{1109415365}{1919384363}\,b \bigg)^2
+\frac{178655579224727}{15976384360300}\bigg(1 \bigg)^2
+\frac{191380511436}{24951996719}\bigg(b \bigg)^2
+\frac{56147093713}{7496706880}\bigg(b^{2}c \bigg)^2.
\end{dmath*}

\begin{dmath*}
h_2(b,c)=\frac{32 b^{10} c^{10}}{15}+\frac{128 b^{10} c^8}{45}+\frac{1936
   b^8 c^{10}}{45}+\frac{13376 b^8 c^8}{225}+\frac{128 b^8
   c^6}{15}+\frac{5104 b^6 c^{10}}{15}+\frac{7424 b^6
   c^8}{15}+\frac{400 b^6 c^6}{3}+\frac{59192 b^4
   c^{10}}{45}+\frac{454432 b^4 c^8}{225}+\frac{32872 b^4
   c^6}{45}+\frac{12374 b^2 c^{10}}{5}+\frac{19896 b^2
   c^8}{5}+\frac{4892 b^2 c^6}{3}+\frac{9009
   c^{10}}{5}+\frac{74484 c^8}{25}+\frac{1}{45} \left(576 b^6
   c^4+3920 b^4 c^4-480 b^4 c^2+5376 b^2 c^4-3210 b^2 c^2+360
   b^2+53838 c^6-4500 c^4-5535 c^2+1080\right).
\end{dmath*}
Using `\texttt{SumsOfSquares}' in \textbf{Macaulay2}, we get
\begin{dmath*}
576 b^6 c^4+3920 b^4 c^4-480 b^4 c^2+5376 b^2 c^4-3210 b^2 c^2+360
   b^2+53838 c^6-4500 c^4-5535 c^2+1080 = 53838 \bigg(\frac{24865}{3876336}\,b^{2}c+c^{3}-\frac{26099}{107676}\,c \bigg)^2
+ 21599 \bigg(-\frac{5209}{194391}\,b^{2}c^{2}+c^{2}-\frac{18043}{86396} \bigg)^2
+\frac{210343}{36} \bigg(\frac{5751}{60098}\,b^{3}c^{2}+b\,c^{2}-\frac{538218}{2734459}\,b \bigg)^2
+\frac{19498324709}{6998076} \bigg(b^{2}c^{2}-\frac{88115791563}{506956442434}\bigg)^2
+\frac{251207361}{480784} \bigg( b^{3}c^{2}-\frac{717899012}{1814275385}\,b\bigg)^2
+\frac{69666947}{215352} \bigg(-\frac{15907018825}{32604131196}\,b^{2}c+c \bigg)^2
+\frac{15599}{52} \bigg(b\,c \bigg)^2
+\frac{11349643727532583}{152587333997280} \bigg(b^{2}c \bigg)^2
+\frac{1418170065555879}{26361735006568} \bigg(1 \bigg)^2
+\frac{42786226455192}{825495300175} \bigg(b \bigg)^2.
\end{dmath*}

\begin{dmath*}
h_3(b,c)=\frac{128 b^{14} c^{14}}{315}+\frac{256 b^{14} c^{12}}{525}+\frac{20672
   b^{12} c^{14}}{1575}+\frac{71552 b^{12} c^{12}}{3675}+\frac{2816
   b^{12} c^{10}}{525}+\frac{11296 b^{10} c^{14}}{63}+\frac{1172672
   b^{10} c^{12}}{3675}+\frac{16544 b^{10} c^{10}}{105}+\frac{2368 b^{10}
   c^8}{105}+\frac{2112752 b^8 c^{14}}{1575}+\frac{3416288 b^8
   c^{12}}{1225}+\frac{988688 b^8 c^{10}}{525}+\frac{252416 b^8
   c^8}{525}+\frac{1152 b^8 c^6}{35}+\frac{373432 b^6
   c^{14}}{63}+\frac{17231216 b^6 c^{12}}{1225}+\frac{1227568 b^6
   c^{10}}{105}+\frac{426784 b^6 c^8}{105}+\frac{18136 b^6
   c^6}{35}+\frac{24322148 b^4 c^{14}}{1575}+\frac{150657352 b^4
   c^{12}}{3675}+\frac{2972504 b^4 c^{10}}{75}+\frac{8826112 b^4
   c^8}{525}+\frac{100708 b^4 c^6}{35}+\frac{766942 b^2
   c^{14}}{35}+\frac{78310956 b^2 c^{12}}{1225}+\frac{2425358 b^2
   c^{10}}{35}+\frac{3529636 b^2 c^8}{105}+\frac{230266 b^2
   c^6}{35}+\frac{325611 c^{14}}{25}+\frac{10067706
   c^{12}}{245}+\frac{8501427 c^{10}}{175}+\frac{4494432 c^8}{175}+\frac{1}{35} \left(784 b^6 c^4+4984 b^4 c^4-560 b^4 c^2+4844
   b^2 c^4-3430 b^2 c^2+420 b^2+177111 c^6-10626 c^4-5495
   c^2+1260\right).
\end{dmath*}
Using `\texttt{SumsOfSquares}' in \textbf{Macaulay2}, we get
\begin{dmath*}
784 b^6 c^4+4984 b^4 c^4-560 b^4 c^2+4844 b^2 c^4-3430 b^2 c^2+420 b^2+177111 c^6-10626 c^4-5495 c^2+1260 = 177111\bigg(\frac{343}{1159272}\,b^{2}c+c^{3}-\frac{470689}{2361480}\,c \bigg)^2
+\frac{1199547}{20}\bigg(-\frac{1090975}{10795923}\,b^{2}c^{2}+c^{2}-\frac{285517}{2399094} \bigg)^2
+\frac{607001}{36}\bigg(-\frac{14247}{6070010}\,b^{3}c^{2}+b\,c^{2}-\frac{1390563}{15782026}\,b \bigg)^2
+\frac{2162205376529}{485816535}\bigg(b^{2}c^{2}-\frac{12112877256075}{56217339789754} \bigg)^2
+\frac{54929220719}{31486400}\bigg(-\frac{4099630565195}{19280156472369}\,b^{2}c+c \bigg)^2
+\frac{190332960599}{242800400}\bigg(b^{3}c^{2}-\frac{1314759474005}{2474328487787}\,b \bigg)^2
+\frac{40527}{52}\bigg(b\,c \bigg)^2
+\frac{5960322065628821313}{29233016690672080}\bigg(1 \bigg)^2
+\frac{703832116675509863}{3759630512111955}\bigg(b^{2}c \bigg)^2
+\frac{4359683320181945}{64332540682462}\bigg(b \bigg)^2.
\end{dmath*}

\begin{dmath*}
h_4(b,c)=\frac{512 b^{18} c^{18}}{14175}+\frac{4096 b^{18}
   c^{16}}{99225}+\frac{7936 b^{16} c^{18}}{3969}+\frac{2725888 b^{16}
   c^{16}}{893025}+\frac{94208 b^{16} c^{14}}{99225}+\frac{4844032 b^{14}
   c^{18}}{99225}+\frac{83216384 b^{14} c^{16}}{893025}+\frac{574976
   b^{14} c^{14}}{11025}+\frac{8192 b^{14} c^{12}}{945}+\frac{68441344
   b^{12} c^{18}}{99225}+\frac{281571328 b^{12}
   c^{16}}{178605}+\frac{8088832 b^{12} c^{14}}{6615}+\frac{4169216
   b^{12} c^{12}}{11025}+\frac{1024 b^{12} c^{10}}{27}+\frac{616465216
   b^{10} c^{18}}{99225}+\frac{14723216384 b^{10}
   c^{16}}{893025}+\frac{529583744 b^{10} c^{14}}{33075}+\frac{231241216
   b^{10} c^{12}}{33075}+\frac{751808 b^{10} c^{10}}{567}+\frac{15872
   b^{10} c^8}{189}+\frac{733857056 b^8 c^{18}}{19845}+\frac{99414094592
   b^8 c^{16}}{893025}+\frac{1413614528 b^8
   c^{14}}{11025}+\frac{334630784 b^8 c^{12}}{4725}+\frac{270950752 b^8
   c^{10}}{14175}+\frac{1192064 b^8 c^8}{525}+\frac{2816 b^8
   c^6}{35}+\frac{14425055648 b^6 c^{18}}{99225}+\frac{435356600192 b^6
   c^{16}}{893025}+\frac{21284903264 b^6 c^{14}}{33075}+\frac{519067904
   b^6 c^{12}}{1225}+\frac{685399136 b^6 c^{10}}{4725}+\frac{115082624
   b^6 c^8}{4725}+\frac{171616 b^6 c^6}{105}+\frac{12032192752 b^4
   c^{18}}{33075}+\frac{1195774968832 b^4
   c^{16}}{893025}+\frac{196261589392 b^4
   c^{14}}{99225}+\frac{49457683168 b^4 c^{12}}{33075}+\frac{8642687696
   b^4 c^{10}}{14175}+\frac{28971328 b^4 c^8}{225}+\frac{59888 b^4
   c^6}{5}+\frac{1158550822 b^2 c^{18}}{2205}+\frac{624204543248 b^2
   c^{16}}{297675}+\frac{22520600344 b^2 c^{14}}{6615}+\frac{19041305312
   b^2 c^{12}}{6615}+\frac{2719211668 b^2 c^{10}}{2025}+\frac{1577198528
   b^2 c^8}{4725}+\frac{1309688 b^2 c^6}{35}+\frac{16338751
   c^{18}}{49}+\frac{28338551384 c^{16}}{19845}+\frac{1107876484
   c^{14}}{441}+\frac{25678259464 c^{12}}{11025}+\frac{5707836706
   c^{10}}{4725}+\frac{532380728 c^8}{1575}+\frac{1}{315} \left(10752
   b^6 c^4+99680 b^4 c^4-6720 b^4 c^2+228256 b^2 c^4-51030 b^2 c^2+5040
   b^2+13164492 c^6-16408 c^4-99435 c^2+20160\right).
\end{dmath*}
Using `\texttt{SumsOfSquares}' in \textbf{Macaulay2}, we get
\begin{dmath*}
10752 b^6 c^4+99680 b^4 c^4-6720 b^4 c^2+228256 b^2 c^4-51030 b^2 c^2+5040 b^2+13164492 c^6-16408 c^4-99435 c^2+20160 = 13164492\bigg(\frac{871111}{473921712}\,b^{2}c+c^{3}-\frac{20957139}{175526560}\,c \bigg)^2
+\frac{62543257}{20}\bigg(-\frac{96138395}{1125778626}\,b^{2}c^{2}+c^{2}-\frac{7591495}{125086514} \bigg)^2
+\frac{25702673}{36}\bigg(-\frac{2527101}{128513365}\,b^{3}c^{2}+b\,c^{2}-\frac{26211375}{668269498}\,b \bigg)^2
+\frac{42535575011281679}{405280305360}\bigg(b^{2}c^{2}-\frac{92045449553337135}{552962475146661827} \bigg)^2
+\frac{649273640274437}{7021062400}\bigg(-\frac{665810223932870}{8440557323567681}\,b^{2}c+c \bigg)^2
+\frac{724643}{26}\bigg(b\,c \bigg)^2
+\frac{26925931846911}{2570267300}\bigg(b^{3}c^{2}-\frac{180424110237905}{350037114009843}\,b \bigg)^2
+\frac{82433649033959276332395}{14377024353813207502}\bigg(1 \bigg)^2
+\frac{1257596798250638235887177}{533274411703006085580}\bigg(b^{2}c \bigg)^2
+\frac{84337814929612787557}{72807719714047344}\bigg(b \bigg)^2.
\end{dmath*}

\begin{dmath*}
h_5(b,c)=\frac{2048 c^{22} b^{22}}{1091475}+\frac{4096 c^{20}
   b^{22}}{1964655}+\frac{31744 c^{22} b^{20}}{200475}+\frac{354304
   c^{20} b^{20}}{1440747}+\frac{53248 c^{18}
   b^{20}}{654885}+\frac{11869696 c^{22} b^{18}}{1964655}+\frac{28832768
   c^{20} b^{18}}{2401245}+\frac{4709888 c^{18}
   b^{18}}{654885}+\frac{31744 c^{16} b^{18}}{24255}+\frac{3336448 c^{22}
   b^{16}}{24255}+\frac{796876288 c^{20} b^{16}}{2401245}+\frac{109120768
   c^{18} b^{16}}{392931}+\frac{186889216 c^{16}
   b^{16}}{1964655}+\frac{806912 c^{14} b^{16}}{72765}+\frac{6798532352
   c^{22} b^{14}}{3274425}+\frac{126837833216 c^{20}
   b^{14}}{21611205}+\frac{12176539136 c^{18}
   b^{14}}{1964655}+\frac{5950041088 c^{16}
   b^{14}}{1964655}+\frac{49040128 c^{14} b^{14}}{72765}+\frac{15872
   c^{12} b^{14}}{297}+\frac{23805268352 c^{22}
   b^{12}}{1091475}+\frac{1521934275328 c^{20}
   b^{12}}{21611205}+\frac{174559620352 c^{18}
   b^{12}}{1964655}+\frac{108932344832 c^{16}
   b^{12}}{1964655}+\frac{775136128 c^{14} b^{12}}{43659}+\frac{5957888
   c^{12} b^{12}}{2205}+\frac{101888 c^{10} b^{12}}{693}+\frac{1315203392
   c^{22} b^{10}}{8085}+\frac{1413228931712 c^{20}
   b^{10}}{2401245}+\frac{1683844485056 c^{18}
   b^{10}}{1964655}+\frac{1264240990336 c^{16}
   b^{10}}{1964655}+\frac{3849765056 c^{14}
   b^{10}}{14553}+\frac{12690841216 c^{12} b^{10}}{218295}+\frac{38016448
   c^{10} b^{10}}{6237}+\frac{6016 c^8 b^{10}}{27}+\frac{1692144108832
   c^{22} b^8}{1964655}+\frac{10617754825024 c^{20}
   b^8}{3087315}+\frac{2217715657952 c^{18}
   b^8}{392931}+\frac{17901245056 c^{16} b^8}{3645}+\frac{531886300384
   c^{14} b^8}{218295}+\frac{4296463040 c^{12} b^8}{6237}+\frac{651507424
   c^{10} b^8}{6237}+\frac{461312 c^8 b^8}{63}+\frac{3328 c^6
   b^8}{21}+\frac{31157070814312 c^{22}
   b^6}{9823275}+\frac{298407089852528 c^{20}
   b^6}{21611205}+\frac{4480724112544 c^{18}
   b^6}{178605}+\frac{48440168876992 c^{16}
   b^6}{1964655}+\frac{3099620236336 c^{14}
   b^6}{218295}+\frac{70563714464 c^{12} b^6}{14553}+\frac{280871456
   c^{10} b^6}{297}+\frac{18037312 c^8 b^6}{189}+\frac{27784 c^6
   b^6}{7}+\frac{25327698773092 c^{22}
   b^4}{3274425}+\frac{785433486144952 c^{20}
   b^4}{21611205}+\frac{141656838089504 c^{18}
   b^4}{1964655}+\frac{154288192539328 c^{16}
   b^4}{1964655}+\frac{67615137704 c^{14} b^4}{1323}+\frac{4421343351248
   c^{12} b^4}{218295}+\frac{29711229488 c^{10} b^4}{6237}+\frac{12907648
   c^8 b^4}{21}+\frac{2280668 c^6 b^4}{63}+\frac{2453547830422 c^{22}
   b^2}{218295}+\frac{81423678540412 c^{20}
   b^2}{1440747}+\frac{15858321568934 c^{18}
   b^2}{130977}+\frac{18863764678564 c^{16}
   b^2}{130977}+\frac{4533671397004 c^{14}
   b^2}{43659}+\frac{10140634746104 c^{12} b^2}{218295}+\frac{7169154292
   c^{10} b^2}{567}+\frac{368840488 c^8 b^2}{189}+\frac{8985182 c^6
   b^2}{63}+\frac{1392896479 c^{22}}{189}+\frac{18906653425790
   c^{20}}{480249}+\frac{3940199652805 c^{18}}{43659}+\frac{5059215383300
   c^{16}}{43659}+\frac{1326722429110 c^{14}}{14553}+\frac{656332750372
   c^{12}}{14553}+\frac{28692391130 c^{10}}{2079}+\frac{51226240
   c^8}{21}+\frac{1}{63} \left(3024 c^4 b^6+37576 c^4 b^4-1680 c^2
   b^4+135100 c^4 b^2-15246 c^2 b^2+1260 b^2+12895025 c^6+107170
   c^4-34923 c^2+6300\right).
\end{dmath*}
Using `\texttt{SumsOfSquares}' in \textbf{Macaulay2}, we get
\begin{dmath*}
3024 c^4 b^6+37576 c^4 b^4-1680 c^2 b^4+135100 c^4 b^2-15246 c^2 b^2+1260 b^2+12895025 c^6+107170 c^4-34923 c^2+6300 = 12895025\bigg(\frac{9158}{12895025}\,b^{2}c+c^{3}-\frac{3751207}{46891000}\,c \bigg)^2
+\frac{43406677}{20}\bigg(-\frac{2898425}{43406677}\,b^{2}c^{2}+c^{2}-\frac{1705331}{43406677} \bigg)^2
+\frac{813253}{2}\bigg(-\frac{51669}{3253012}\,b^{3}c^{2}+b\,c^{2}-\frac{69487}{3020654}\,b \bigg)^2
+\frac{99568634438661}{1875640000}\bigg(-\frac{2717579463920}{33189544812887}\,b^{2}c+c \bigg)^2
+\frac{1771700819399}{43406677}\bigg(b^{2}c^{2}-\frac{6844893001259}{46064221304374} \bigg)^2
+\frac{186727}{13}\bigg(b\,c \bigg)^2
+\frac{76027180743}{26024096}\bigg(b^{3}c^{2}-\frac{481698502754}{988353349659}\,b \bigg)^2
+\frac{3067336004076266963}{1497087192392155}\bigg(1 \bigg)^2
+\frac{4822899680844187461}{9492209816485682}\bigg(b^{2}c \bigg)^2
+\frac{4508374085874757}{12848593545567}\bigg(b \bigg)^2.
\end{dmath*}

\begin{footnotesize}\begin{dmath*}
h_6(b,c)=\frac{8192 c^{26} b^{26}}{127702575}+\frac{32768 c^{24}
   b^{26}}{468242775}+\frac{10760192 c^{26}
   b^{24}}{1404728325}+\frac{73252864 c^{24}
   b^{24}}{6087156075}+\frac{1933312 c^{22}
   b^{24}}{468242775}+\frac{590606336 c^{26}
   b^{22}}{1404728325}+\frac{248741888 c^{24}
   b^{22}}{289864575}+\frac{251260928 c^{22}
   b^{22}}{468242775}+\frac{4407296 c^{20}
   b^{22}}{42567525}+\frac{3945715712 c^{26}
   b^{20}}{280945665}+\frac{213729304576 c^{24}
   b^{20}}{6087156075}+\frac{14447581184 c^{22}
   b^{20}}{468242775}+\frac{1759670272 c^{20}
   b^{20}}{156080925}+\frac{581632 c^{18}
   b^{20}}{405405}+\frac{149143711232 c^{26}
   b^{18}}{468242775}+\frac{634940704768 c^{24}
   b^{18}}{676350675}+\frac{32664988672 c^{22}
   b^{18}}{31216185}+\frac{17120804864 c^{20}
   b^{18}}{31216185}+\frac{42243584 c^{18} b^{18}}{315315}+\frac{108544
   c^{16} b^{18}}{9009}+\frac{346414848256 c^{26}
   b^{16}}{66891825}+\frac{35540563938304 c^{24}
   b^{16}}{2029052025}+\frac{332493773312 c^{22}
   b^{16}}{14189175}+\frac{492099224576 c^{20}
   b^{16}}{31216185}+\frac{15809266432 c^{18}
   b^{16}}{2837835}+\frac{549275648 c^{16} b^{16}}{567567}+\frac{2854912
   c^{14} b^{16}}{45045}+\frac{87228518486528 c^{26}
   b^{14}}{1404728325}+\frac{1438014162632704 c^{24}
   b^{14}}{6087156075}+\frac{1465996959232 c^{22}
   b^{14}}{4002075}+\frac{46494374803456 c^{20}
   b^{14}}{156080925}+\frac{42837903872 c^{18}
   b^{14}}{315315}+\frac{10812772352 c^{16}
   b^{14}}{315315}+\frac{117697024 c^{14} b^{14}}{27027}+\frac{149504
   c^{12} b^{14}}{715}+\frac{31236824525056 c^{26}
   b^{12}}{56189133}+\frac{14257521908077568 c^{24}
   b^{12}}{6087156075}+\frac{1923969255517952 c^{22}
   b^{12}}{468242775}+\frac{55196549321216 c^{20}
   b^{12}}{14189175}+\frac{6107609131264 c^{18}
   b^{12}}{2837835}+\frac{1995434665984 c^{16}
   b^{12}}{2837835}+\frac{123057842944 c^{14}
   b^{12}}{945945}+\frac{3828476416 c^{12} b^{12}}{315315}+\frac{1457152
   c^{10} b^{12}}{3465}+\frac{5219339619628192 c^{26}
   b^{10}}{1404728325}+\frac{11581321577165696 c^{24}
   b^{10}}{676350675}+\frac{240480639853184 c^{22}
   b^{10}}{7203735}+\frac{671544878167040 c^{20}
   b^{10}}{18729711}+\frac{7298447665216 c^{18}
   b^{10}}{315315}+\frac{17565080320 c^{16}
   b^{10}}{1911}+\frac{2085803000192 c^{14}
   b^{10}}{945945}+\frac{94694944768 c^{12}
   b^{10}}{315315}+\frac{70219552 c^{10} b^{10}}{3465}+\frac{30592 c^8
   b^{10}}{63}+\frac{1978323914357264 c^{26}
   b^8}{108056025}+\frac{556119878319132224 c^{24}
   b^8}{6087156075}+\frac{91412424181657792 c^{22}
   b^8}{468242775}+\frac{2428744315669888 c^{20}
   b^8}{10405395}+\frac{161932407186016 c^{18}
   b^8}{945945}+\frac{45073399934720 c^{16}
   b^8}{567567}+\frac{21927261204928 c^{14}
   b^8}{945945}+\frac{86134081408 c^{12} b^8}{21021}+\frac{279871952
   c^{10} b^8}{693}+\frac{395072 c^8 b^8}{21}+\frac{1920 c^6
   b^8}{7}+\frac{3363043725212368 c^{26}
   b^6}{52026975}+\frac{18017954680636544 c^{24}
   b^6}{52026975}+\frac{41738823811154992 c^{22}
   b^6}{52026975}+\frac{164252072011062592 c^{20}
   b^6}{156080925}+\frac{271117905758368 c^{18}
   b^6}{315315}+\frac{142978544588288 c^{16}
   b^6}{315315}+\frac{29267159829280 c^{14}
   b^6}{189189}+\frac{10494171683456 c^{12} b^6}{315315}+\frac{2949785008
   c^{10} b^6}{693}+\frac{18196352 c^8 b^6}{63}+\frac{57328 c^6
   b^6}{7}+\frac{690479779159336 c^{26}
   b^4}{4459455}+\frac{54371705319629344 c^{24}
   b^4}{61486425}+\frac{114494495345469496 c^{22}
   b^4}{52026975}+\frac{23269671261024112 c^{20}
   b^4}{7432425}+\frac{295079874070736 c^{18}
   b^4}{105105}+\frac{520193894301824 c^{16}
   b^4}{315315}+\frac{605240132759056 c^{14}
   b^4}{945945}+\frac{50714422114784 c^{12}
   b^4}{315315}+\frac{86651069368 c^{10} b^4}{3465}+\frac{46056544 c^8
   b^4}{21}+\frac{625944 c^6 b^4}{7}+\frac{3192312960134 c^{26}
   b^2}{14157}+\frac{2929503725483144 c^{24}
   b^2}{2147145}+\frac{4185541698679412 c^{22}
   b^2}{1156155}+\frac{19193256200856064 c^{20}
   b^2}{3468465}+\frac{37709362565806 c^{18}
   b^2}{7007}+\frac{24315798069560 c^{16}
   b^2}{7007}+\frac{157559074357192 c^{14}
   b^2}{105105}+\frac{1363235388608 c^{12} b^2}{3185}+\frac{29878490722
   c^{10} b^2}{385}+\frac{24684424 c^8 b^2}{3}+424156 c^6
   b^2+\frac{164190371975 c^{26}}{1089}+\frac{5361237368884
   c^{24}}{5577}+\frac{89260211464222
   c^{22}}{33033}+\frac{145446654518708
   c^{20}}{33033}+\frac{4604119725105 c^{18}}{1001}+\frac{22528817322120
   c^{16}}{7007}+\frac{2910348074588 c^{14}}{1911}+\frac{3389265520952
   c^{12}}{7007}+\frac{1092754339 c^{10}}{11}+\frac{85356220
   c^8}{7}+\frac{1}{7} \left(448 c^4 b^6+7056 c^4 b^4-224 c^2 b^4+34496
   c^4 b^2-2366 c^2 b^2+168 b^2+5169090 c^6+48636 c^4-6223
   c^2+1008\right).
\end{dmath*}\end{footnotesize}
Using `\texttt{SumsOfSquares}' in \textbf{Macaulay2}, we get
\begin{dmath*}
448 c^4 b^6+7056 c^4 b^4-224 c^2 b^4+34496 c^4 b^2-2366 c^2 b^2+168 b^2+5169090 c^6+48636 c^4-6223 c^2+1008 = 5169090\bigg(-\frac{17579}{41352720}\,b^{2}c+c^{3}-\frac{160928}{2584545}\,c \bigg)^2
+692348\bigg(-\frac{79731}{1384696}\,b^{2}c^{2}+c^{2}-\frac{153597}{5538784} \bigg)^2
+\frac{474487}{4}\bigg(-\frac{61767}{4744870}\,b^{3}c^{2}+b\,c^{2}-\frac{200809}{12336662}\,b \bigg)^2
+\frac{125461294753}{10338180}\bigg(-\frac{218508674731}{3261993663578}\,b^{2}c+c \bigg)^2
+\frac{108683246871}{13846960}\bigg(b^{2}c^{2}-\frac{735071362015}{5651528837292} \bigg)^2
+\frac{139289}{52}\bigg(b\,c \bigg)^2
+\frac{81212908111}{189794800}\bigg(b^{3}c^{2}-\frac{555976793995}{1055767805443}\,b \bigg)^2
+\frac{402956990227706381}{1175517998156736}\bigg(1 \bigg)^2
+\frac{1642106627418391347}{13569893640484480}\bigg(b^{2}c \bigg)^2
+\frac{245774985506386}{13724981470759}\bigg(b \bigg)^2.
\end{dmath*}

\begin{scriptsize}\begin{dmath*}
h_7(b,c)=\frac{32768 c^{30} b^{30}}{21070924875}+\frac{65536 c^{28}
   b^{30}}{39131717625}+\frac{68337664 c^{30}
   b^{28}}{273922023375}+\frac{231964672 c^{28}
   b^{28}}{586975764375}+\frac{5439488 c^{26}
   b^{28}}{39131717625}+\frac{242819072 c^{30}
   b^{26}}{13043905875}+\frac{22728261632 c^{28}
   b^{26}}{586975764375}+\frac{108470272 c^{26}
   b^{26}}{4347968625}+\frac{5029888 c^{24}
   b^{26}}{1003377375}+\frac{33536528384 c^{30}
   b^{24}}{39131717625}+\frac{430012915712 c^{28}
   b^{24}}{195658588125}+\frac{646418432 c^{26}
   b^{24}}{323402625}+\frac{1106870272 c^{24}
   b^{24}}{1449322875}+\frac{20676608 c^{22}
   b^{24}}{200675475}+\frac{152180918272 c^{30}
   b^{22}}{5590245375}+\frac{4411522617344 c^{28}
   b^{22}}{53361433125}+\frac{35678670848 c^{26}
   b^{22}}{372683025}+\frac{686989582336 c^{24}
   b^{22}}{13043905875}+\frac{329676800 c^{22}
   b^{22}}{24081057}+\frac{4870144 c^{20}
   b^{22}}{3648645}+\frac{1177566524416 c^{30}
   b^{20}}{1863415125}+\frac{1297201921439744 c^{28}
   b^{20}}{586975764375}+\frac{120590942124032 c^{26}
   b^{20}}{39131717625}+\frac{3164546449408 c^{24}
   b^{20}}{1449322875}+\frac{496066573312 c^{22}
   b^{20}}{602026425}+\frac{128628736 c^{20}
   b^{20}}{825825}+\frac{9900032 c^{18}
   b^{20}}{868725}+\frac{3033593405456896 c^{30}
   b^{18}}{273922023375}+\frac{8547949053494272 c^{28}
   b^{18}}{195658588125}+\frac{923760806332928 c^{26}
   b^{18}}{13043905875}+\frac{793678304392192 c^{24}
   b^{18}}{13043905875}+\frac{17900065622528 c^{22}
   b^{18}}{602026425}+\frac{1644496237568 c^{20}
   b^{18}}{200675475}+\frac{460692992 c^{18}
   b^{18}}{394875}+\frac{18885632 c^{16}
   b^{18}}{289575}+\frac{40866338200811264 c^{30}
   b^{16}}{273922023375}+\frac{29504266738674176 c^{28}
   b^{16}}{45151981875}+\frac{47023328208361216 c^{26}
   b^{16}}{39131717625}+\frac{5241665987442688 c^{24}
   b^{16}}{4347968625}+\frac{15985056058624 c^{22}
   b^{16}}{22297275}+\frac{2452215712256 c^{20}
   b^{16}}{9555975}+\frac{231990900992 c^{18}
   b^{16}}{4343625}+\frac{8173356032 c^{16}
   b^{16}}{1403325}+\frac{72306688 c^{14}
   b^{16}}{289575}+\frac{6772078994287744 c^{30}
   b^{14}}{4347968625}+\frac{4390359989064600832 c^{28}
   b^{14}}{586975764375}+\frac{13327300113976832 c^{26}
   b^{14}}{869593725}+\frac{228173159255886848 c^{24}
   b^{14}}{13043905875}+\frac{349050803730176 c^{22}
   b^{14}}{28667925}+\frac{3906839916032 c^{20}
   b^{14}}{735075}+\frac{6252159817216 c^{18}
   b^{14}}{4343625}+\frac{4194160620544 c^{16}
   b^{14}}{18243225}+\frac{206661248 c^{14} b^{14}}{10725}+\frac{12118784
   c^{12} b^{14}}{19305}+\frac{493086967257621952 c^{30}
   b^{12}}{39131717625}+\frac{12849033838793635712 c^{28}
   b^{12}}{195658588125}+\frac{5786157093705459712 c^{26}
   b^{12}}{39131717625}+\frac{91015857306318848 c^{24}
   b^{12}}{483107625}+\frac{29987049724226944 c^{22}
   b^{12}}{200675475}+\frac{48708220722944 c^{20}
   b^{12}}{637065}+\frac{9984511641344 c^{18}
   b^{12}}{394875}+\frac{2742910942208 c^{16}
   b^{12}}{521235}+\frac{62455173824 c^{14} b^{12}}{96525}+\frac{53144704
   c^{12} b^{12}}{1287}+\frac{1478912 c^{10}
   b^{12}}{1485}+\frac{437984906344618208 c^{30}
   b^{10}}{5590245375}+\frac{257911544235495186752 c^{28}
   b^{10}}{586975764375}+\frac{182138381290570528 c^{26}
   b^{10}}{169401375}+\frac{1787221027036076224 c^{24}
   b^{10}}{1185809625}+\frac{14622490648806592 c^{22}
   b^{10}}{10945935}+\frac{6258176637986432 c^{20}
   b^{10}}{8027019}+\frac{1313921511295936 c^{18}
   b^{10}}{4343625}+\frac{40108717234304 c^{16}
   b^{10}}{521235}+\frac{1191295347488 c^{14}
   b^{10}}{96525}+\frac{22338232768 c^{12} b^{10}}{19305}+\frac{244108192
   c^{10} b^{10}}{4455}+\frac{125504 c^8
   b^{10}}{135}+\frac{685089599910084784 c^{30}
   b^8}{1863415125}+\frac{1293046639384693665376 c^{28}
   b^8}{586975764375}+\frac{17475723595240975024 c^{26}
   b^8}{3010132125}+\frac{12828101842262460032 c^{24}
   b^8}{1449322875}+\frac{472692113777782048 c^{22}
   b^8}{54729675}+\frac{125670015478511936 c^{20}
   b^8}{22297275}+\frac{3613520248080416 c^{18}
   b^8}{1447875}+\frac{1942340622708352 c^{16}
   b^8}{2606175}+\frac{14112701617264 c^{14}
   b^8}{96525}+\frac{345486271904 c^{12} b^8}{19305}+\frac{1110551824
   c^{10} b^8}{891}+\frac{1037824 c^8 b^8}{25}+\frac{2176 c^6
   b^8}{5}+\frac{38357652830898352024 c^{30}
   b^6}{30435780375}+\frac{1568896744214797657648 c^{28}
   b^6}{195658588125}+\frac{58930304512834223408 c^{26}
   b^6}{2608781175}+\frac{37232233827217872032 c^{24}
   b^6}{1003377375}+\frac{23716853781261738904 c^{22}
   b^6}{602026425}+\frac{5681394540067554544 c^{20}
   b^6}{200675475}+\frac{6767389044544096 c^{18}
   b^6}{482625}+\frac{12458458699837888 c^{16}
   b^6}{2606175}+\frac{106297404992888 c^{14}
   b^6}{96525}+\frac{1061306887088 c^{12} b^6}{6435}+\frac{7424792752
   c^{10} b^6}{495}+\frac{495362848 c^8 b^6}{675}+\frac{226648 c^6
   b^6}{15}+\frac{6044717066970154732 c^{30}
   b^4}{2029052025}+\frac{1305734628921976084072 c^{28}
   b^4}{65219529375}+\frac{260462198364442750552 c^{26}
   b^4}{4347968625}+\frac{457436751807728788544 c^{24}
   b^4}{4347968625}+\frac{24154389256946167852 c^{22}
   b^4}{200675475}+\frac{899441721348971752 c^{20}
   b^4}{9555975}+\frac{223013428893905216 c^{18}
   b^4}{4343625}+\frac{71435302008678208 c^{16}
   b^4}{3648645}+\frac{1491052459583948 c^{14}
   b^4}{289575}+\frac{17479572907528 c^{12}
   b^4}{19305}+\frac{448607662984 c^{10} b^4}{4455}+\frac{1443895168 c^8
   b^4}{225}+\frac{2882692 c^6 b^4}{15}+\frac{1864510962275218 c^{30}
   b^2}{429429}+\frac{423442590899827348 c^{28}
   b^2}{13803075}+\frac{267667308866139454 c^{26}
   b^2}{2760615}+\frac{1497908249959569716 c^{24}
   b^2}{8281845}+\frac{2114042008345258174 c^{22}
   b^2}{9555975}+\frac{593790439445694364 c^{20}
   b^2}{3185325}+\frac{14562325156401034 c^{18}
   b^2}{131625}+\frac{282509860040384828 c^{16}
   b^2}{6081075}+\frac{1318793421067142 c^{14}
   b^2}{96525}+\frac{3528101498444 c^{12} b^2}{1287}+\frac{1592710244282
   c^{10} b^2}{4455}+\frac{18733155076 c^8 b^2}{675}+\frac{3194002 c^6
   b^2}{3}+\frac{60139207926675 c^{30}}{20449}+\frac{444912639060106
   c^{28}}{20449}+\frac{1478290911554065
   c^{26}}{20449}+\frac{8737529975507480
   c^{24}}{61347}+\frac{238105585232221
   c^{22}}{1287}+\frac{2359065426708998
   c^{20}}{14157}+\frac{186978113252539
   c^{18}}{1755}+\frac{563604806416076
   c^{16}}{11583}+\frac{20269646811239 c^{14}}{1287}+\frac{7572883528966
   c^{12}}{2145}+\frac{778679579839 c^{10}}{1485}+\frac{10684326464
   c^8}{225}+\frac{1}{45} \left(3696 c^4 b^6+70952 c^4 b^4-1680 c^2
   b^4+438004 c^4 b^2-20250 c^2 b^2+1260 b^2+97752669 c^6+847658
   c^4-60105 c^2+8820\right).
\end{dmath*}\end{scriptsize}
Using `\texttt{SumsOfSquares}' in \textbf{Macaulay2}, we get
\begin{dmath*}
3696 c^4 b^6+70952 c^4 b^4-1680 c^2 b^4+438004 c^4 b^2-20250 c^2 b^2+1260 b^2+97752669 c^6+847658 c^4-60105 c^2+8820 = 97752669\bigg(-\frac{23981}{97752669}\,b^{2}c+c^{3}-\frac{49373207}{977526690}\,c \bigg)^2
+\frac{53611497}{5}\bigg(-\frac{11286025}{214445988}\,b^{2}c^{2}+c^{2}-\frac{1006625}{47654664} \bigg)^2
+\frac{3229137}{2}\bigg(-\frac{94043}{7175860}\,b^{3}c^{2}+b\,c^{2}-\frac{758857}{55971708}\,b \bigg)^2
+\frac{350689158239069}{2443816725}\bigg(-\frac{793505948202085}{18235836228431588}\,b^{2}c+c \bigg)^2
+\frac{358437370956251}{4288919760}\bigg(b^{2}c^{2}-\frac{1009088034241725}{9319371644862526} \bigg)^2
+\frac{1448441}{52}\bigg(b\,c \bigg)^2
+\frac{981282369759}{287034400}\bigg(b^{3}c^{2}-\frac{6562235249365}{12756670806867}\,b \bigg)^2
+\frac{2961857827865059555905}{969214651065702704}\bigg(1 \bigg)^2
+\frac{934901241025865976851}{948263483878442576}\bigg(b^{2}c \bigg)^2
+\frac{19419295884942205}{331673440978542}\bigg(b \bigg)^2.
\end{dmath*}

\section{Proof of HFRI (\ref{Aug18a}) for the case $m_2\ge 8$} \label{sec4} \setcounter{equation}{0}

We will show that  (\ref{Aug18a}) holds for $m_2\ge 8$. To this end, we split the unit interval into two subintervals according to the bound
\begin{eqnarray}\label{BBB}
B=\frac{2.75}{m_2m_3}.
\end{eqnarray}
In the following two subsections, we use different techniques to prove that
\begin{eqnarray}\label{Aug26A}
\frac{F\left(-m_2,-m_3;\frac{1}{2};z\right)}{F\left(-m_2,-m_3;\frac{3}{2};z\right)}>\frac{1}{H(z)}
\end{eqnarray}
for the case $z\le \frac{2.75}{m_2m_3}$ and the case $z> \frac{2.75}{m_2m_3}$ separately.

Here we would like to point out that although we choose to use the bound (\ref{BBB}) corresponding to the truncation number ``4" (see (\ref{Aug26CC}) below), we may also use the following bound corresponding to the truncation number ``5":
\begin{eqnarray}\label{BBB2}
B=\frac{3.8}{m_2m_3}.
\end{eqnarray}
The advantage of using the bound  (\ref{BBB2}) is that we can apply the same method of this section to establish the inequality (\ref{Aug26A}) for all $m_2\ge 6$. The price is that the Appendix would become much longer. To  shorten the length of this paper, we have decided to use the truncation number ``4" and the bound  (\ref{BBB}).

\subsection{The case that $z\le \frac{2.75}{m_2m_3}$ }\label{41}

First, we present a useful lemma.
\begin{lem} For $m_2,m_3\in\mathbb{N}$, we have
\begin{eqnarray}\label{Aug199}
H(z)>
 \begin{cases}
 \frac{1}{2},\ \  &{\rm if}\ \frac{1}{r^2}< z\le\frac{1}{r},\\
 \frac{1}{7},\ \ &{\rm if}\ \frac{1}{r}< z\le\frac{2.75}{m_2m_3}.
   \end{cases}
\end{eqnarray}
\end{lem}

\noindent{\bf Proof.}\ \ For  $\frac{1}{r^2}< z\le\frac{1}{r}$, we have
\begin{eqnarray*}
&&H(z)>\frac{1}{2}\\
&\Leftrightarrow&4[(m_3-m_2)(rz-1)]^2+4(r-1)^3z>[(r^2z-1)+2(m_2+m_3+1)(1-rz)]^2\\
&\Leftrightarrow&4(r-1)^3z-(r^2z-1)^2-4(r-1)(rz-1)^2>4(m_2+m_3+1)(1-rz)(r^2z-1)\\
&\Leftrightarrow&4(r-1)(1-z)(r^2z-1)-(r^2z-1)^2>4(m_2+m_3+1)(1-rz)(r^2z-1)\\
&\Leftrightarrow&4(r-1)(1-z)-(r^2z-1)>4(m_2+m_3+1)(1-rz)\\
&\Leftarrow&3r-7>4(m_2+m_3+1)(1-rz)\\
&\Leftarrow&3r-7>4(m_2+m_3+1).
\end{eqnarray*}
Then,
\begin{eqnarray*}
H(z)>\frac{1}{2},\ \ \ \ \frac{1}{r^2}< z\le\frac{1}{r}.
\end{eqnarray*}

For  $\frac{1}{r}< z\le\frac{2.75}{m_2m_3}$, we have
\begin{eqnarray*}
&&H(z)>\frac{1}{7}\\
&\Leftrightarrow&49[(m_3-m_2)(rz-1)]^2+49(r-1)^3z>[(r^2z-1)-7(m_2+m_3+1)(rz-1)]^2\\
&\Leftrightarrow&49(r-1)^3z-(r^2z-1)^2-49(r-1)(rz-1)^2>-14(m_2+m_3+1)(rz-1)(r^2z-1)\\
&\Leftrightarrow&49(r-1)(1-z)(r^2z-1)-(r^2z-1)^2>-14(m_2+m_3+1)(rz-1)(r^2z-1)\\
&\Leftrightarrow&49(r-1)(1-z)-(r^2z-1)>-14(m_2+m_3+1)(rz-1)\\
&\Leftrightarrow&49r-r^2z+49z+[14(m_2+m_3+1)-49]rz>14(m_2+m_3+1)+48\\
&\Leftarrow&14.5r>14(m_2+m_3+1)+48.
\end{eqnarray*}
Then,
\begin{eqnarray*}
H(z)>\frac{1}{7},\ \ \ \ \frac{1}{r}< z\le\frac{2.75}{m_2m_3}.
\end{eqnarray*}\hfill\fbox

By (\ref{Aug27HH}), we obtain that for $\frac{1}{r^2}< z<1$,
\begin{eqnarray}\label{Aug26CC}
&&\frac{F\left(-m_2,-m_3;\frac{1}{2};z\right)}{F\left(-m_2,-m_3;\frac{3}{2};z\right)}-\frac{1}{H(z)}>0\nonumber\\
&\Leftrightarrow&\ \ \left[(m_2+m_3+1)(rz-1)+\sqrt{[(m_3-m_2)(rz-1)]^2+(2m_2+1)^3(2m_3+1)^3z}\right]\nonumber\\
&&\ \ \cdot\sum_{j=0}^{m_2}\frac{2^{2j}z^j}{(m_2-j)!(m_3-j)!(2j)!}\nonumber\\
&&>(r^2z-1)\sum_{j=0}^{m_2}\frac{2^{2j}z^j}{(m_2-j)!(m_3-j)!(2j+1)!}\nonumber\\
&\Leftarrow&\ \ \ \ \ \ \ \ \ \ \left[(m_2+m_3+1)(rz-1)+\sqrt{[(m_3-m_2)(rz-1)]^2+(2m_2+1)^3(2m_3+1)^3z}\right]\nonumber\\
&&\ \ \ \ \ \ \ \ \ \ \cdot\sum_{j=0}^{4}\frac{2^{2j}z^j}{(m_2-j)!(m_3-j)!(2j)!}\nonumber\\
&&\ \ \ \ \ \ \ >(r^2z-1)\sum_{j=0}^{4}\frac{2^{2j}z^j}{(m_2-j)!(m_3-j)!(2j+1)!},\nonumber\\
&&{\rm and}\ \ H(z)\sum_{j=5}^{m_2}\frac{2^{2j}z^j}{(m_2-j)!(m_3-j)!(2j)!}>\sum_{j=5}^{m_2}\frac{2^{2j}z^j}{(m_2-j)!(m_3-j)!(2j+1)!}.
\end{eqnarray}

Define
$$
u=m_2m_3z.
$$
Similar to (\ref{Aug22g}), we can show that
\begin{eqnarray}\label{Aug20aa}
&&\ \ \left[(m_2+m_3+1)(rz-1)+\sqrt{[(m_3-m_2)(rz-1)]^2+(2m_2+1)^3(2m_3+1)^3z}\right]\nonumber\\
&&\ \ \cdot\left[\sum_{j=0}^{4}\frac{2^{2j}z^j}{(m_2-j)!(m_3-j)!(2j)!}\right]\nonumber\\
&&>(r^2z-1)\left[\sum_{j=0}^{4}\frac{2^{2j}z^j}{(m_2-j)!(m_3-j)!(2j+1)!}\right]\nonumber\\
&\Leftarrow&\ \ (2m_2+1)(2m_3+1)\left(1-\frac{u}{m_2m_3}\right)\left[1+\sum_{j=1}^{4}\frac{2^{2j}u^j(m_2-1)!(m_3-1)!}{(m_2m_3)^{j-1}(m_2-j)!(m_3-j)!(2j)!}\right]^2\nonumber\\
&&\ \ +2(m_2+m_3+1)\left\{\frac{[(2m_2+1)(2m_3+1)+1]u}{m_2m_3}-1\right\}\nonumber\\
&&\ \ \ \ \cdot\left[1+\sum_{j=1}^{4}\frac{2^{2j}u^j(m_2-1)!(m_3-1)!}{(m_2m_3)^{j-1}(m_2-j)!(m_3-j)!(2j)!}\right]\nonumber\\
&&\ \ \ \ \cdot\left[1+\sum_{j=1}^{4}\frac{2^{2j}u^j(m_2-1)!(m_3-1)!}{(m_2m_3)^{j-1}(m_2-j)!(m_3-j)!(2j+1)!}\right]\nonumber\\
&&>\left\{\frac{[(2m_2+1)(2m_3+1)+1]^2u}{m_2m_3}-1\right\}\nonumber\\
&&\ \ \cdot\left[1+\sum_{j=1}^{4}\frac{2^{2j}u^j(m_2-1)!(m_3-1)!}{(m_2m_3)^{j-1}(m_2-j)!(m_3-j)!(2j+1)!}\right]^2.
\end{eqnarray}
Obviously,  the second inequality of (\ref{Aug20aa}) is a direct consequence of the following inequality: for $u\in(0,2.75)$ and $x_2\ge8$, $x_3\ge x_2$,
\begin{eqnarray}\label{Aug29l}
&&f(x_2,x_3,u)\nonumber\\
&:=&(2x_2+1)(2x_3+1)\left({x_2x_3}-{u}\right)\left[17!(x_2x_3)^3+\sum_{j=1}^{4}\frac{2^{2j}u^j(x_2x_3)^{4-j}(x_2-1)!(x_3-1)!17!}{(x_2-j)!(x_3-j)!(2j)!}\right]^2\nonumber\\
&&+2(x_2+x_3+1)\left\{{[(2x_2+1)(2x_3+1)+1]u}-{x_2x_3}\right\}\nonumber\\
&&\ \ \cdot\left[17!(x_2x_3)^3+\sum_{j=1}^{4}\frac{2^{2j}u^j(x_2x_3)^{4-j}(x_2-1)!(x_3-1)!17!}{(x_2-j)!(x_3-j)!(2j)!}\right]\nonumber\\
&&\ \ \cdot\left[17!(x_2x_3)^3+\sum_{j=1}^{4}\frac{2^{2j}u^j(x_2x_3)^{4-j}(x_2-1)!(x_3-1)!17!}{(x_2-j)!(x_3-j)!(2j+1)!}\right]\nonumber\\
&&-\left\{{[(2x_2+1)(2x_3+1)+1]^2u}-{x_2x_3}\right\}\nonumber\\
&&\ \ \cdot\left[17!(x_2x_3)^3+\sum_{j=1}^{4}\frac{2^{2j}u^j(x_2x_3)^{4-j}(x_2-1)!(x_3-1)!17!}{(x_2-j)!(x_3-j)!(2j+1)!}\right]^2\nonumber\\
&>&0.
\end{eqnarray}

By using the transformations
$$
u=\frac{2.75c^2}{1+c^2},\ \ x_2=a^2+8,\ \ x_3=b^2+8,
$$
we define
$$
g(a,b,c):=(1+c^2)^{9}\cdot f(x_2,x_3,u),\ \ \ \ a,b,c \in \mathbb{R}.
$$
Then, (\ref{Aug29l}) holds for any $u\in(0,2.75)$ and $x_2\ge8$, $x_3\ge x_2$ if $g$ is positive on $\mathbb{R}^3$.
By virtue of \textbf{Mathematica}, we obtain the expansion of $g$. See the Appendix.
Note that all terms in the  expansion are positive.  Hence (\ref{Aug29l})  holds. Therefore, the first inequality of  (\ref{Aug20aa}) holds, which together with (\ref{Aug199}) and (\ref{Aug26CC}) implies that
\begin{eqnarray*}
\frac{F\left(-m_2,-m_3;\frac{1}{2};z\right)}{F\left(-m_2,-m_3;\frac{3}{2};z\right)}>\frac{1}{H(z)},\ \ \ \ t<z\le \frac{2.75}{m_2m_3}.
\end{eqnarray*}

\subsection{The case that $z> \frac{2.75}{m_2m_3}$ }

In this subsection, we show that for $m_2\ge 8$,
\begin{eqnarray*}
\frac{F\left(-m_2,-m_3;\frac{1}{2};z\right)}{F\left(-m_2,-m_3;\frac{3}{2};z\right)}>\frac{1}{H(z)},\ \ \ \ \frac{2.75}{m_2m_3}<z< 1.
\end{eqnarray*}
Note that
\begin{eqnarray*}
&&\frac{F\left(-m_2,-m_3;\frac{1}{2};z\right)}{F\left(-m_2,-m_3;\frac{3}{2};z\right)}>\frac{1}{H(z)}\nonumber\\
&\Leftrightarrow&F\left(-m_2,-m_3;\frac{3}{2};z\right)-H(z)F\left(-m_2,-m_3;\frac{1}{2};z\right)<0\nonumber\\
&\Leftrightarrow&\frac{F\left(-m_2-1,-m_3;\frac{1}{2};z\right)-(1-z)F\left(-m_2,-m_3;\frac{1}{2};z\right)}{(2m_3+1)z}-H(z)F\left(-m_2,-m_3;\frac{1}{2};z\right)<0\nonumber\\
&\Leftrightarrow&F\left(-m_2-1,-m_3;\frac{1}{2};z\right)-[(1-z)+(2m_3+1)zH(z)]F\left(-m_2,-m_3;\frac{1}{2};z\right)<0.
\end{eqnarray*}
For $0<z\le 1$, define
\begin{eqnarray}\label{Aug15aaa}
G(z):=F\left(-m_2-1,-m_3;\frac{1}{2};z\right)-[(1-z)+(2m_3+1)zH(z)]F\left(-m_2,-m_3;\frac{1}{2};z\right).
\end{eqnarray}
We will show that
\begin{eqnarray}\label{Aug15a}
G(z)<0,\ \ \ \ \frac{2.75}{m_2m_3}< z<1.
\end{eqnarray}

By (\ref{KKJJGG}) and  (\ref{HGF}), we get
{\small\begin{eqnarray}\label{HGFD}
&&H(1)\nonumber\\
&=&\lim_{z\rightarrow 1}\frac{1-z}{-\beta+\sqrt{\beta^2-(1-z)\gamma}}\nonumber\\
&=&\lim_{z\rightarrow 1}\frac{1-z}{(m_2+m_3+1)\left[
\frac{1-z}{(2m_2+1)(2m_3+1)}-z\right]+\sqrt{(m_3-m_2)^2\left[
\frac{1-z}{(2m_2+1)(2m_3+1)}-z\right]^2+(2m_2+1)(2m_3+1)z}}\nonumber\\
&=&\lim_{z\rightarrow 1}\{(1-z)\nonumber\\
&&\left.\ \ \cdot\frac{(m_2+m_3+1)\left[
\frac{1-z}{(2m_2+1)(2m_3+1)}-z\right]-\sqrt{(m_3-m_2)^2\left[
\frac{1-z}{(2m_2+1)(2m_3+1)}-z\right]^2+(2m_2+1)(2m_3+1)z}}{(2m_2+1)(2m_3+1)\left[
\frac{1-z}{(2m_2+1)(2m_3+1)}-z\right]^2-(2m_2+1)(2m_3+1)z}\right\}\nonumber\\
&=&\lim_{z\rightarrow 1}\frac{(m_2+m_3+1)\left[
\frac{1-z}{(2m_2+1)(2m_3+1)}-z\right]-\sqrt{(m_3-m_2)^2\left[
\frac{1-z}{(2m_2+1)(2m_3+1)}-z\right]^2+(2m_2+1)(2m_3+1)z}}{\frac{1-z}{(2m_2+1)(2m_3+1)}-2z-(2m_2+1)(2m_3+1)z}\nonumber\\
&=&\frac{2(m_2+m_3+1)}{2+(2m_2+1)(2m_3+1)}.
\end{eqnarray}}

\noindent Then,
\begin{eqnarray*}
G(1)&=&F\left(-m_2-1,-m_3;\frac{1}{2};1\right)-(2m_3+1)\cdot H(1)\cdot F\left(-m_2,-m_3;\frac{1}{2};1\right)\\
&=&\frac{\Gamma(\frac{1}{2})\Gamma(m_2+m_3+\frac{3}{2})}{\Gamma(m_2+\frac{3}{2})\Gamma(m_3+\frac{1}{2})}-\frac{2(2m_3+1)(m_2+m_3+1)}{2+(2m_2+1)(2m_3+1)}\cdot\frac{\Gamma(\frac{1}{2})\Gamma(m_2+m_3+\frac{1}{2})}{\Gamma(m_2+\frac{1}{2})\Gamma(m_3+\frac{1}{2})}\\
&=&-\frac{(2m_2-1)(2m_3-1)-2}{(2m_2+1)[2+(2m_2+1)(2m_3+1)]}\cdot\frac{\Gamma(\frac{1}{2})\Gamma(m_2+m_3+\frac{1}{2})}{\Gamma(m_2+\frac{1}{2})\Gamma(m_3+\frac{1}{2})}\\
&<&0.
\end{eqnarray*}
Hence, to prove (\ref{Aug15a}), we need only show that for $\frac{2.75}{m_2m_3}< z<1$,
$$
G'(z)=0\Rightarrow G(z)<0.
$$

By (\ref{Gauss3}), we get
\begin{eqnarray}\label{Aug15j}
&&G'(z)=0\nonumber\\
&\Leftrightarrow&\left\{F\left(-m_2-1,-m_3;\frac{1}{2};z\right)-[(1-z)+(2m_3+1)zH(z)]F\left(-m_2,-m_3;\frac{1}{2};z\right)\right\}'=0\nonumber\\
&\Leftrightarrow&0=\frac{(m_2+1)\left[F\left(-m_2-1,-m_3;\frac{1}{2};z\right)-F\left(-m_2,-m_3;\frac{1}{2};z\right)\right]}{z}\nonumber\\
&&\ \ \ \ \ -[(1-z)+(2m_3+1)zH(z)]'F\left(-m_2,-m_3;\frac{1}{2};z\right)\nonumber\\
&&\ \ \ \ \ -[(1-z)+(2m_3+1)zH(z)]\frac{m_2\left[F\left(-m_2,-m_3;\frac{1}{2};z\right)-F\left(-m_2+1,-m_3;\frac{1}{2};z\right)\right]}{z}\nonumber\\
&\Leftrightarrow&0=\frac{(m_2+1)\left[F\left(-m_2-1,-m_3;\frac{1}{2},z\right)-F\left(-m_2,-m_3;\frac{1}{2};z\right)\right]}{z}\nonumber\\
&&\ \ \ \ \ -[(1-z)+(2m_3+1)zH(z)]'F\left(-m_2,-m_3;\frac{1}{2};z\right)\nonumber\\
&&\ \ \ \ \ -[(1-z)+(2m_3+1)zH(z)]\nonumber\\
&&\ \ \ \ \ \cdot\frac{(2m_2+1)F\left(-m_2-1,-m_3;\frac{1}{2};z\right)-[(2m_2+1)+2m_3z]F\left(-m_2,-m_3;\frac{1}{2};z\right)}{2z(1-z)}\nonumber\\
&\Leftrightarrow&\ \ \ \{2(m_2+1)(1-z)-(2m_2+1)[(1-z)+(2m_3+1)zH(z)]\}F\left(-m_2-1,-m_3;\frac{1}{2};z\right)\nonumber\\
&&=\{2(m_2+1)(1-z)+2z(1-z)[(1-z)+(2m_3+1)zH(z)]'\nonumber\\
&&\ \ \ \ -[(2m_2+1)+2m_3z][(1-z)+(2m_3+1)zH(z)]\}F\left(-m_2,-m_3;\frac{1}{2};z\right)\nonumber\\
&\Leftrightarrow&\ \ \ \{(1-z)-(2m_2+1)(2m_3+1)zH(z)\}F\left(-m_2-1,-m_3;\frac{1}{2};z\right)\nonumber\\
&&=\{[1-2(m_3+1)z](1-z)+2(2m_3+1)z^2(1-z)H'(z)\nonumber\\
&&\ \ \ \ +(2m_3+1)z[1-2m_2-2(m_3+1)z]H(z)\}F\left(-m_2,-m_3;\frac{1}{2};z\right).
\end{eqnarray}
Note that  for $\frac{2.75}{m_2m_3}< z<1$,
\begin{eqnarray*}
&&(1-z)-(2m_2+1)(2m_3+1)zH(z)\\
&<&1-\frac{[(2m_2+1)(2m_3+1)]^{5/2}z^{3/2}}{r^2z-1}\\
&<&1-\frac{[(2m_2+1)(2m_3+1)]^{5/2}z^{1/2}}{[(2m_2+1)(2m_3+1)+1]^2}\\
&<&0.
\end{eqnarray*}
Thus, by (\ref{Aug15aaa}) and (\ref{Aug15j}), to complete the proof of (\ref{Aug15a}) it suffices to show that for $\frac{2.75}{m_2m_3}< z<1$,
\begin{eqnarray*}
&&[1-2(m_3+1)z](1-z)+2(2m_3+1)z^2(1-z)H'(z)+(2m_3+1)z[1-2m_2-2(m_3+1)z]H(z)\nonumber\\
&>&[(1-z)-(2m_2+1)(2m_3+1)zH(z)]\cdot[(1-z)+(2m_3+1)zH(z)],
\end{eqnarray*}
i.e.,
\begin{eqnarray}\label{Aug13v}
(1-z)+[2(m_2+m_3+1)z-1]H(z)<(2m_2+1)(2m_3+1)zH^2(z)+2z(1-z)H'(z).
\end{eqnarray}

We have
\begin{eqnarray*}
H'(z)&=&\frac{1}{2(r^2z-1)^2\sqrt{[(m_3-m_2)(rz-1)]^2+(2m_2+1)^3(2m_3+1)^3z}}\\
&&\cdot\left\{2r(r-1)(m_2+m_3+1)\sqrt{[(m_3-m_2)(rz-1)]^2+(2m_2+1)^3(2m_3+1)^3z}\right.\\
&&\ \ -(1+r^2z)(2m_2+1)^3(2m_3+1)^3+2(m_3-m_2)^2r(r-1)(rz-1)\}.
\end{eqnarray*}
Then, for $\frac{2.75}{m_2m_3}< z<1$, we get
\begin{eqnarray*}
&&(\ref{Aug13v})\ {\rm holds}\\
&\Leftrightarrow&\ \ (r-1)z\left\{\frac{(m_2+m_3+1)(rz-1)+\sqrt{[(m_3-m_2)(rz-1)]^2+(r-1)^3z}}{r^2z-1}\right\}^2\nonumber\\
&&\ \ +\frac{z(1-z)}{(r^2z-1)^2\sqrt{[(m_3-m_2)(rz-1)]^2+(r-1)^3z}}\nonumber\\
&&\ \ \ \ \cdot\left\{2r(r-1)(m_2+m_3+1)\sqrt{[(m_3-m_2)(rz-1)]^2+(r-1)^3z}+2(m_3-m_2)^2r(r-1)(rz-1)\right\}\nonumber\\
&&>\frac{(r-1)^3z(1-z)(1+r^2z)}{(r^2z-1)^2\sqrt{[(m_3-m_2)(rz-1)]^2+(r-1)^3z}}
+(1-z)\nonumber\\
&&\ \ +[2(m_2+m_3+1)z-1]\cdot\frac{(m_2+m_3+1)(rz-1)+\sqrt{[(m_3-m_2)(rz-1)]^2+(r-1)^3z}}{r^2z-1}\nonumber\\
&\Leftrightarrow&\ \ \frac{ (r-1)z\{[(r-1)+2(m_3-m_2)^2](rz-1)^2+(r-1)^3z\}}{(r^2z-1)^2}\nonumber\\
&&\ \ +\frac{2(m_2+m_3+1)(r-1)z(rz-1)\sqrt{[(m_3-m_2)(rz-1)]^2+(r-1)^3z}}{(r^2z-1)^2}\nonumber\\
&&\ \ +\frac{2(m_2+m_3+1)r(r-1)z(1-z)}{(r^2z-1)^2}\nonumber\\
&&\ \ +\frac{2(m_3-m_2)^2r(r-1)(rz-1)z(1-z)}{(r^2z-1)^2\sqrt{[(m_3-m_2)(rz-1)]^2+(r-1)^3z}}\nonumber\\
&&>\frac{(r-1)^3z(1-z)(1+r^2z)}{(r^2z-1)^2\sqrt{[(m_3-m_2)(rz-1)]^2+(r-1)^3z}}
+(1-z)\nonumber\\
&&\ \ +\frac{[2(m_2+m_3+1)z-1]\left[(m_2+m_3+1)(rz-1)+\sqrt{[(m_3-m_2)(rz-1)]^2+(r-1)^3z}\right]}{r^2z-1}\nonumber\\
&\Leftrightarrow&\ \ \frac{(r-1)^4z^2+(m_2+m_3+1)[2r(r-1)z(1-z)+(rz-1)(r^2z-1)]}{(r^2z-1)^2}\nonumber\\
&&>\frac{(r-1)z(1-z)\left[(r-1)^2(1+r^2z)-2(m_3-m_2)^2r(rz-1)\right]}{(r^2z-1)^2\sqrt{[(m_3-m_2)(rz-1)]^2+(r-1)^3z}}
+(1-z)\nonumber\\
&&\ \ +\frac{[2(m_2+m_3+1)z(r+rz-2)-(r^2z-1)]\sqrt{[(m_3-m_2)(rz-1)]^2+(r-1)^3z}}{(r^2z-1)^2}\nonumber\\
&&\ \ +\frac{(r-1)z(rz-1)(r^2z+rz+r-3)+2(m_3-m_2)^2z(rz-1)(rz+r-2)}{(r^2z-1)^2}\\
&\Leftrightarrow&\ \ (r-1)^4z^2+(m_2+m_3+1)[2r(r-1)z(1-z)+(rz-1)(r^2z-1)]\nonumber\\
&&>\frac{(r-1)z(1-z)\left[(r-1)^2(1+r^2z)-2(m_3-m_2)^2r(rz-1)\right]}{\sqrt{[(m_3-m_2)(rz-1)]^2+(r-1)^3z}}
+(r^2z-1)^2(1-z)\nonumber\\
&&\ \ +[2(m_2+m_3+1)z(r+rz-2)-(r^2z-1)]\sqrt{[(m_3-m_2)(rz-1)]^2+(r-1)^3z}\nonumber\\
&&\ \ +(r-1)z(rz-1)(r^2z+rz+r-3)+2(m_3-m_2)^2z(rz-1)(rz+r-2)\nonumber
\end{eqnarray*}
\begin{eqnarray}\label{RRRR}
&\Leftrightarrow&\ \ -1 + 4 z - 4 r z + 3 r^2 z + z^2 - 8 r z^2 + 8 r^2 z^2 - 4 r^3 z^2 +
 r^2 z^3\nonumber\\
&&\ \ +(m_2+m_3+1)(1 - 3 r z + r^2 z + 2 r z^2 - 2 r^2 z^2 + r^3 z^2)\nonumber\\
&&\ \ -2(m_3-m_2)^2(2 z - r z - 3 r z^2 + r^2 z^2 + r^2 z^3)\nonumber\\
&&>\frac{1}{\sqrt{[(m_3-m_2)(rz-1)]^2+(r-1)^3z}}\nonumber\\
&&\ \ \cdot \left\{(r-1)z(1-z)\left[(r-1)^2(1+r^2z)-2(m_3-m_2)^2r(rz-1)\right]\right.\nonumber\\
&&\ \ \left.+[2(m_2+m_3+1)z(r+rz-2)-(r^2z-1)]\{[(m_3-m_2)(rz-1)]^2+(r-1)^3z\}\right\}.\nonumber\\
&&
\end{eqnarray}

For  $\frac{2.75}{m_2m_3}< z<1$ , we have
\begin{eqnarray*}
&&-1 + 4 z - 4 r z + 3 r^2 z + z^2 - 8 r z^2 + 8 r^2 z^2 - 4 r^3 z^2 +
 r^2 z^3\nonumber\\
&&+(m_2+m_3+1)(1 - 3 r z + r^2 z + 2 r z^2 - 2 r^2 z^2 + r^3 z^2)\nonumber\\
&&-2(m_3-m_2)^2(2 z - r z - 3 r z^2 + r^2 z^2 + r^2 z^3)\nonumber\\
&>&- 4 r^3 z^2+\sqrt{(m_3-m_2)^2+(r-1)}(- 3 r^2 z^2 + r^3 z^2)-4(m_3-m_2)^2 r^2z^2\\
&\ge&- 4 r^3 z^2+\sqrt{(m_3-m_2)^2+(r-1)}\cdot\frac{287r^3 z^2}{290} -4(m_3-m_2)^2 r^2z^2\\
&>&\left(\frac{1}{2}\sqrt{(m_3-m_2)^2+(r-1)}-4\right)r^3 z^2+\frac{142(m_3-m_2)r^3 z^2}{290}-4(m_3-m_2)^2 r^2z^2\\
&>&0,
\end{eqnarray*}
and, for $z\ge \frac{2.1}{2m_2+1}$, we have
\begin{eqnarray*}
&&(r-1)z(1-z)\left[(r-1)^2(1+r^2z)-2(m_3-m_2)^2r(rz-1)\right]\\
&&+[2(m_2+m_3+1)z(r+rz-2)-(r^2z-1)]\{[(m_3-m_2)(rz-1)]^2+(r-1)^3z\}\\
&<&\left(1-\frac{2.1}{2m_2+1}\right)(r-1)^3z(r^2z+1)+\{4(m_2+m_3+1)rz-(r^2z-1)\}(r-1)^3z\\
&=&\left\{2\left(1-\frac{2.1}{2m_2+1}\right)+4(m_2+m_3+1)rz-\frac{2.1(r^2z-1)}{2m_2+1}\right\}(r-1)^3z\\
&=&\left\{\frac{2-\frac{2.1}{2m_2+1}}{\frac{2.1r^2z}{2m_2+1}}+\frac{4(m_2+m_3+1)rz}{\frac{2.1r^2z}{2m_2+1}}-1\right\}\frac{2.1r^2(r-1)^3z^2}{2m_2+1}\\
&<&\left\{\frac{2}{(2.1)^2(2m_3+1)^2}+\frac{2}{2.1}-1\right\}\frac{2.1r^2(r-1)^3z^2}{2m_2+1}\\
&<&0.
\end{eqnarray*}
Then, by (\ref{RRRR}), we conclude that (\ref{Aug13v}) holds for $\frac{2.1}{2m_2+1}\le z<1$.

Finally, we will show that (\ref{RRRR}) holds for $\frac{2.75}{m_2m_3}< z<\frac{2.1}{2m_2+1}$. To this end, we assume without loss of generality that the right hand side of  (\ref{RRRR}) is positive. Then, for $\frac{2.75}{m_2m_3}< z<\frac{2.1}{2m_2+1}$, we have
\begin{eqnarray*}
&&(\ref{RRRR})\ {\rm holds}\\
&\Leftarrow&\ \ \{-1 + 4 z - 4 r z + 3 r^2 z + z^2 - 8 r z^2 + 8 r^2 z^2 - 4 r^3 z^2 +
 r^2 z^3\\
&&\ \ +\sqrt{(m_3-m_2)^2+(r-1)}(1 - 3 r z + r^2 z + 2 r z^2 - 2 r^2 z^2 + r^3 z^2)\\
&&\ \ -2(m_3-m_2)^2(2 z - r z - 3 r z^2 + r^2 z^2 + r^2 z^3)\}^2\{[(m_3-m_2)(rz-1)]^2+(r-1)^3z\}\\
&&>\left\{(r-1)z(1-z)\left[(r-1)^2(1+r^2z)-2(m_3-m_2)^2r(rz-1)\right]\right.\\
&&\ \ \left.+[2\sqrt{(m_3-m_2)^2+(r-1)}z(r+rz-2)-(r^2z-1)]\{[(m_3-m_2)(rz-1)]^2+(r-1)^3z\}\right\}^2\\
&\Leftrightarrow&\ \ [(m_3-m_2)^2 (r z-1)^2+(
       r-1)^3 z  ]\\
&&\ \ \cdot \left\{-1 + 4 z - 4 r z + 3 r^2 z + z^2 -
    8 r z^2 + 8 r^2 z^2 - 4 r^3 z^2 + r^2 z^3\right.\\
&&\ \ \ \ \ +
    \sqrt{(m_3-m_2)^2+(r-1)
      } [1 + r^3 z^2 + r z (-3 + 2 z) + r^2 (z - 2 z^2)]\\
&&\ \ \ \ \ \left.-
    2 (m_3-m_2)^2 z [2 + r^2 z (1 + z) - r (1 + 3 z)]\right\}^2\\
&&\ \ -\left\{[(m_3-m_2)^2 ( r z-1)^2+(r-1 )^3 z ] [1+
       2 \sqrt{(m_3-m_2)^2+(r-1) } z (r + r z-2  )]\right. \\
&&\ \ \ \ \ \ -(m_3-m_2)^2rz ( r z-1)[2(r-1)(1-z)+r(rz-1)] \\
&&\ \ \ \ \ \ \left.-(r-1)^3r^2z^3+(r-1)^3z(1-z)\right\}^2\\
&&>0.
\end{eqnarray*}
Therefore, the proof of (\ref{RRRR}) for $\frac{2.75}{m_2m_3}< z<\frac{2.1}{2m_2+1}$ is complete by
\begin{eqnarray*}
&& [(m_3-m_2)^2 (r z-1)^2+(
       r-1)^3 z  ]\\
&&\ \ \cdot \left\{-1 + 4 z - 4 r z + 3 r^2 z + z^2 -
    8 r z^2 + 8 r^2 z^2 - 4 r^3 z^2 + r^2 z^3\right.\\
&&\ \ \ \ \ +
    \sqrt{(m_3-m_2)^2+(r-1)
      } [1 + r^3 z^2 + r z (-3 + 2 z) + r^2 (z - 2 z^2)]\\
&&\ \ \ \ \ \left.-
    2 (m_3-m_2)^2 z [2 + r^2 z (1 + z) - r (1 + 3 z)]\right\}^2\\
&&\ \ -\left\{[(m_3-m_2)^2 ( r z-1)^2+(r-1 )^3 z ] [1+
       2 \sqrt{(m_3-m_2)^2+(r-1) } z (r + r z-2  )]\right. \\
&&\ \ \ \ \ \ -(m_3-m_2)^2rz ( r z-1)[2(r-1)(1-z)+r(rz-1)] \\
&&\ \ \ \ \ \ \left.-(r-1)^3r^2z^3+(r-1)^3z(1-z)\right\}^2\\
&>&(       r-1)^3 z \left\{ - 4 r^3 z^2 +
    \sqrt{(m_3-m_2)^2+(r-1)
      } (r^3 z^2 -3rz)-2 (m_3-m_2)^2  r^2 z^2(1+z)\right\}^2\\
&&-\left\{(r-1 )^3 z (2-z-r^2z^2)+
       2 \sqrt{(m_3-m_2)^2+(r-1) }(r-1 )^3 rz^2 (1 + z)\right\}^2\\
&>&( r-1)^3 r^6z^5 [(m_3-m_2)^2+(r-1)
      ]\left\{1-\frac{4}{(r-1)^{1/2}} - \frac{3}{r^2 z}- \frac{2(m_3-m_2)(1+z)}{r}\right\}^2\\
&&-4[(m_3-m_2)^2+(r-1) ](r-1)^6r^2z^4 (1 + z)^2\\
&>& [(m_3-m_2)^2+(r-1)](r-1 )^6 r^2 z ^4\\
&&\cdot\left[\left\{1-\frac{4}{(r-1)^{1/2}} - \frac{3}{r^2 z}- \frac{1+\frac{2.1}{2m_2+1}}{2m_2+1}\right\}^2rz-4(1 + z)^2\right]\\
&>& [(m_3-m_2)^2+(r-1)](r-1 )^6 r^2 z ^4\\
&&\cdot\left[\left\{1-\frac{4}{17} - \frac{3}{(17^2)(4)(2.75)}- \frac{1+\frac{2.1}{17}}{17}\right\}^2( 4)( 2.75)-4\left(1 + \frac{2.1}{17}\right)^2\right]\\
&>&0.3[(m_3-m_2)^2+(r-1)](r-1 )^6 r^2 z ^4\\
&>&0.
\end{eqnarray*}

\section{Concluding remarks}\label{sec5}\setcounter{equation}{0}

{\bf Remark 1}\ \ In this paper, the exponents $m_2$ and $m_3$ are assumed to be natural numbers. However,  Theorems \ref{thmm1m2m3} and \ref{thm2} can be extended to the case that  $m_2$ and $m_3$ are positive real numbers.

Let $y_2,y_3 \in (0,\infty)$. We consider a centered Gaussian random vector $(X_1,X_2,X_3)$. Let $X_1 = X_2+aX_3$ for some $a \in \mathbb{R}$. Define $x=E[X_2 X_3]$. Assume without loss of generality that $E[X_2^2]=E[X_3^2]=1$. We have
\begin{equation}\label{1d}
    E[|X_2|^{y_2}]=\frac{2^{y_2/2}\Gamma(\frac{y_2+1}{2})}{\sqrt{\pi}},\ \ \ \ E[|X_3|^{y_3}]=\frac{2^{y_3/2}\Gamma(\frac{y_3+1}{2})}{\sqrt{\pi}}.
\end{equation}
By Nabeya \cite{Nabeya1951}, we get
\begin{eqnarray}\label{2deven}
&&E[|X_2|^{y_2} |X_3|^{y_3+2}] = (y_3+1)\cdot\frac{2^{(y_2+y_3)/2}\Gamma(\frac{y_2+1}{2})\Gamma(\frac{y_3+1}{2})}{\pi}F\left(-\frac{y_3}{2}-1,-\frac{y_2}{2};\frac{1}{2};x^2\right),\nonumber\\
    &&E[|X_2|^{y_2+2} |X_3|^{y_3}] = (y_2+1)\cdot\frac{2^{(y_2+y_3)/2}\Gamma(\frac{y_2+1}{2})\Gamma(\frac{y_3+1}{2})}{\pi}F\left(-\frac{y_2}{2}-1,-\frac{y_3}{2};\frac{1}{2};x^2\right).\ \ \ \
\end{eqnarray}
Let $p(x_2,x_3)$ be the probability density function of $(X_2,X_3)$. By Kamat \cite{Kamat1958}, we get
\begin{eqnarray}\label{2dodd}
&&E[|X_2|^{y_2}X_2 |X_3|^{y_3}X_3] \nonumber\\
&=& \int_{-\infty}^{\infty}\int_{-\infty}^{\infty}|x_2|^{y_2} x_2 |x_3|^{y_3} x_3 p(x_2,x_3) dx_2 dx_3 \nonumber\\
    &=&\left[\int_{0}^{\infty}\int_{0}^{\infty}+\int_{-\infty}^{0}\int_{-\infty}^{0}-\int_{0}^{\infty}\int_{-\infty}^{0}- \int_{-\infty}^{0}\int_{0}^{\infty}\right]|x_2|^{y_2+1} |x_3|^{y_3+1} p(x_2,x_3) dx_2 dx_3\nonumber\\
    &=&2\cdot\frac{2^{(y_2+y_3-2)/2}}{\pi}\Bigg[\Gamma\left(\frac{y_2+2}{2}\right)\Gamma\left(\frac{y_3+2}{2}\right)F\left(-\frac{y_2+1}{2},-\frac{y_3+1}{2};\frac{1}{2};x^2\right)\nonumber\\
    &&\ \ \ \ \ \ \ \ \ \ \ \ \ \ \ \ \ \ \ \ +2x\Gamma\left(\frac{y_2+3}{2}\right)\Gamma\left(\frac{y_3+3}{2}\right)F\left(-\frac{y_2}{2},-\frac{y_3}{2};\frac{3}{2};x^2\right)\Bigg]\nonumber\\
    &&-2\cdot\frac{2^{(y_2+y_3-2)/2}}{\pi}\Bigg[\Gamma\left(\frac{y_2+2}{2}\right)\Gamma\left(\frac{y_3+2}{2}\right)F\left(-\frac{y_2+1}{2},-\frac{y_3+1}{2};\frac{1}{2};x^2\right)\nonumber\\
    &&\ \ \ \ \ \ \ \ \ \ \ \ \ \ \ \ \ \ \ \ \ -2x\Gamma\left(\frac{y_2+3}{2}\right)\Gamma\left(\frac{y_3+3}{2}\right)F\left(-\frac{y_2}{2},-\frac{y_3}{2};\frac{3}{2};x^2\right)\Bigg]\nonumber\\
    &&=x(y_2+1)(y_3+1)\cdot\frac{2^{(y_2+y_3)/2}\Gamma\left(\frac{y_2+1}{2}\right)\Gamma\left(\frac{y_3+1}{2}\
\right)}{\pi}F\left(-\frac{y_2}{2},-\frac{y_3}{2};\frac{3}{2};x^2\right).
\end{eqnarray}
Then,
\begin{eqnarray*}
&&\frac{E[(X_2+aX_3)^{2}|X_2|^{y_2} |X_3|^{y_3} ]} {E[X_1^{2}]E[|X_2|^{y_2}] E[|X_3|^{y_3}] }\nonumber\\
&=&\frac{E[a^2|X_2|^{y_2} |X_3|^{y_3+2}] + |X_2|^{y_2+2} |X_3|^{y_3} + 2a|X_2|^{y_2}X_2 |X_3|^{y_3}X_3]} {(a^2+1+2ax)E[|X_2|^{y_2}] E[|X_3|^{y_3}]}\nonumber\\
&=&\Bigg[a^2\left(y_3+1\right)F\left(-\frac{y_3}{2}-1,-\frac{y_2}{2};\frac{1}{2};x^2\right)+\left(y_2+1\right)F\left(-\frac{y_3}{2},-\frac{y_2}{2}-1;\frac{1}{2};x^2\right)\nonumber\\
&&+2ax\left(y_3+1\right)\left(y_2+1\right) F\left(-\frac{y_3}{2},-\frac{y_2}{2};\frac{3}{2};x^2\right)\Bigg]\cdot\frac{1}{a^2+1+2ax}
\end{eqnarray*}
by (\ref{1d})--(\ref{2dodd}). Hence, by a modified version of \cite[Lemma 2.1]{RusSun}, verifying the corresponding GPIs is equivalent to showing that for any $a\in\mathbb{R}$ and $x\in[-1,1]$,
\begin{eqnarray}\label{gen1}
&&a^2\left(y_3+1\right)F\left(-\frac{y_3}{2}-1,-\frac{y_2}{2};\frac{1}{2};x^2\right)+\left(y_2+1\right)F\left(-\frac{y_3}{2},-\frac{y_2}{2}-1;\frac{1}{2};x^2\right)\nonumber\\
&&+2ax\left(y_3+1\right)\left(y_2+1\right) F\left(-\frac{y_3}{2},-\frac{y_2}{2};\frac{3}{2};x^2\right)\nonumber\\
&>&a^2+1+2ax.
\end{eqnarray}

Note that the inequality (\ref{gen1}) is exactly the same as the inequality (\ref{gen2}) except that integer-valued exponents $m_2,m_3$ are replaced with real-valued exponents $\frac{y_2}{2},\frac{y_3}{2}$. The arguments used in Sections \ref{sec2} and \ref{sec4} can be applied without any change. Additionally, as explained in the beginning of Section \ref{sec4}, we can cover the case  that
$\frac{y_2}{2},\frac{y_3}{2}\ge 6$ by using the truncation number ``5" and the bound $B=\frac{3.8}{\frac{y_2}{2}\cdot\frac{y_3}{2}}$. Therefore, we have the following propositions.
\begin{pro}\label{thmm1m2m3U}
Let $y_2,y_3 \in [12,\infty)$. For any centered Gaussian random vector $(X_1,X_2,X_3)$,
\begin{eqnarray}\label{GPIjkl}
E[X_1^{2} |X_2|^{y_2} |X_3|^{y_3}]\ge E[X_1^{2}] E[|X_2|^{y_2}] E[|X_3|^{y_3}].
\end{eqnarray}
The equality sign holds if and only if $X_1,X_2,X_3$ are independent.
\end{pro}

Define
$$
r_{y_2,y_3}=(y_2+1)(y_3+1)+1,\ \ \ \ t_{y_2,y_3}=\frac{1}{r_{y_2,y_3}+\left(1+\frac{1}{y_2}\right)\left(1+\frac{1}{y_3}\right)},
$$
and  for $\frac{1}{r^2_{y_2,y_3}}<z\le 1$,
\begin{eqnarray*}
H_{y_2,y_3}(z)=\frac{\frac{(y_2+y_3+2)(r_{y_2,y_3}z-1)}{2}+\sqrt{\frac{[(y_3-y_2)(r_{y_2,y_3}z-1)]^2}{4}+(y_2+1)^3(y_3+1)^3z}}{r_{y_2,y_3}^2z-1}.
\end{eqnarray*}

\begin{pro}\label{thm2V}
Let $(X_2,X_3)$ be a centered Gaussian random vector. If $y_2,y_3\in [12,\infty)$, then
\begin{eqnarray}\label{MRIjkl}
&&\frac{\left|E[ |X_2|^{y_2}X_2|X_3|^{y_3}X_3]\right|}{(y_2+1)(y_3+1)E[ |X_2|^{y_2}|X_3|^{y_3}]}\nonumber\\
&\le&
 \begin{cases}
 |{\rm Cov}(X_2,X_3)|,\ \  &{\rm if}\ |{\rm Corr}(X_2,X_3)|\le \sqrt{t_{y_2,y_3}},\nonumber\\
H_{y_2,y_3}([{\rm Corr}(X_2,X_3)]^2)\cdot|{\rm Cov}(X_2,X_3)|,\ \ &{\rm if}\ \sqrt{t_{y_2,y_3}}<|{\rm Corr}(X_2,X_3)|.
   \end{cases}\nonumber\\
&&
\end{eqnarray}
The equality sign holds if and only if $X_2$ and $X_3$ are independent.
\end{pro}

Through more delicate analysis, the interval $[12,\infty)$ used in Propositions \ref{thmm1m2m3U} and \ref{thm2V} can be enlarged. However, Proposition \ref{thm2V} does not hold for all $y_2,y_3\in(0,\infty)$. For example, \textbf{Mathematica} has shown that  the inequality  (\ref{MRIjkl}) does not hold when $y_2=4$ and $y_3\in[4,4.58]$, or when $y_2=2$ and $y_3\in[2,8.15]$. Therefore, the validity of the GPI (\ref{GPIjkl}) for the most general case that $y_2,y_3\in (0,\infty)$ still remains open, although using  \textbf{Mathematica} to check the inequality (\ref{gen1}), it appears to be true.
\vskip 0.5cm

\noindent {\bf Remark 2}\ \ It is natural to ask if the MRI method developed in this paper can be adapted to solve the following 3D-GPI.

\begin{con}\label{conmm}
Let $m_1,m_2,m_3\in\mathbb{N}$. For any centered Gaussian random vector $(X_1,X_2,X_3)$,
$$
E[X_1^{2m_1} X_2^{2m_2}X_3^{2m_3}]\ge E[X_1^{2m_1}] E[X_2^{2m_2}] E[X_3^{2m_3}].
$$
The equality sign holds if and only if $X_1,X_2,X_3$ are independent.
\end{con}

By a modified version of \cite[Lemma 2.1]{RusSun}, we may assume without loss of generality that $X_1= X_2+aX_3$ for some $a \in \mathbb{R}$ and $E[X_2^2]=E[X_3^2]=1$. Define $x=E[X_2 X_3]$. Similar to (\ref{moment}), we consider the moment ratio:
\begin{eqnarray*}
&&\frac{E[(X_2+aX_3)^{2m_1}X_2^{2m_2}X_3^{2m_3} ]} {E[X_1^{2m_1}]E[X_2^{2m_2}] E[X_3^{2m_3}] }\nonumber\\
&=&\frac{\sum_{k=0}^{2m_1}\binom{2m_1}{k}a^kE[X_2^{2m_2+2m_1-k}X_3^{2m_3+k} ]} {(a^2+1+2ax)^{m_1}(2m_1-1)!!(2m_2-1)!!(2m_3-1)!!}\nonumber\\
&=&\Bigg[\sum_{i=0}^{m_1}\binom{2m_1}{2i}a^{2i}(2m_3+2i-1)\cdots(2m_3+1)\cdot(2m_2+2m_1-2i-1)\cdots(2m_2+1)\nonumber\\
&&\ \ \cdot F\left(-m_3-i,-m_2-m_1+i;\frac{1}{2};x^2\right)\nonumber\\
&&+x\sum_{j=1}^{m_1}\binom{2m_1}{2j-1}a^{2j-1}(2m_3+2j-1)\cdots(2m_3+1)\cdot(2m_2+2m_1-2j+1)\cdots(2m_2+1)\nonumber\\
&&\ \ \cdot F\left(-m_3-j+1,-m_2-m_1+j;\frac{3}{2};x^2\right)\Bigg]\cdot\frac{1}{(a^2+1+2ax)^{m_1}(2m_1-1)!!}.
\end{eqnarray*}
The proof of Conjecture \ref{conmm} is complete if we can show that
\begin{eqnarray}\label{gen3}
&&\sum_{i=0}^{m_1}\binom{2m_1}{2i}a^{2i}(2m_3+2i-1)\cdots(2m_3+1)\cdot(2m_2+2m_1-2i-1)\cdots(2m_2+1)\nonumber\\
&&\ \ \cdot  F\left(-m_3-i,-m_2-m_1+i;\frac{1}{2};x^2\right)\nonumber\\
&&+x\sum_{j=1}^{m_1}\binom{2m_1}{2j-1}a^{2j-1}(2m_3+2j-1)\cdots(2m_3+1)\cdot(2m_2+2m_1-2j+1)\cdots(2m_2+1)\nonumber\\
&&\ \ \cdot F\left(-m_3-j+1,-m_2-m_1+j;\frac{3}{2};x^2\right)\nonumber\\
&>&(a^2+1+2ax)^{m_1}(2m_1-1)!!.
\end{eqnarray}
Obviously, the inequality (\ref{gen3}) extends the inequality (\ref{gen2}).  Further, by using the  relations (\ref{Gauss3}), we may use only two hypergeometric functions, $F\left(-m_2,-m_3;\frac{1}{2};x^2\right)$ and $F\left(-m_2,-m_3;\frac{3}{2};x^2\right)$, to simplify (\ref{gen3}).
We hope the MRI method introduced in this paper can be improved so as to  prove (\ref{gen3}) with all three exponents unbounded.

\section{Appendix: Expansion of the function $g$ in \S \ref{41}}\label{sec6}\setcounter{equation}{0}

We have the following SOS expression of the function $g$ for any $a,b,c \in \mathbb{R}$:
\begin{verbatim}
 g(a,b,c)=148260632637820250986905600 + 148260632637820250986905600 a^2 +
 64864026779046359806771200 a^4 + 16216006694761589951692800 a^6 +
 2533751046056498429952000 a^8 + 253375104605649842995200 a^10 +
 15835944037853115187200 a^12 + 565569429923325542400 a^14 +
 8837022342551961600 a^16 + 148260632637820250986905600 b^2 +
 148260632637820250986905600 a^2 b^2 +
 64864026779046359806771200 a^4 b^2 +
 16216006694761589951692800 a^6 b^2 +
 2533751046056498429952000 a^8 b^2 +
 253375104605649842995200 a^10 b^2 +
 15835944037853115187200 a^12 b^2 + 565569429923325542400 a^14 b^2 +
 8837022342551961600 a^16 b^2 + 64864026779046359806771200 b^4 +
 64864026779046359806771200 a^2 b^4 +
 28378011715832782415462400 a^4 b^4 +
 7094502928958195603865600 a^6 b^4 +
 1108516082649718063104000 a^8 b^4 +
 110851608264971806310400 a^10 b^4 +
 6928225516560737894400 a^12 b^4 + 247436625591454924800 a^14 b^4 +
 3866197274866483200 a^16 b^4 + 16216006694761589951692800 b^6 +
 16216006694761589951692800 a^2 b^6 +
 7094502928958195603865600 a^4 b^6 +
 1773625732239548900966400 a^6 b^6 +
 277129020662429515776000 a^8 b^6 +
 27712902066242951577600 a^10 b^6 + 1732056379140184473600 a^12 b^6 +
 61859156397863731200 a^14 b^6 + 966549318716620800 a^16 b^6 +
 2533751046056498429952000 b^8 + 2533751046056498429952000 a^2 b^8 +
 1108516082649718063104000 a^4 b^8 +
 277129020662429515776000 a^6 b^8 + 43301409478504611840000 a^8 b^8 +
 4330140947850461184000 a^10 b^8 + 270633809240653824000 a^12 b^8 +
 9665493187166208000 a^14 b^8 + 151023331049472000 a^16 b^8 +
 253375104605649842995200 b^10 + 253375104605649842995200 a^2 b^10 +
 110851608264971806310400 a^4 b^10 +
 27712902066242951577600 a^6 b^10 + 4330140947850461184000 a^8 b^10 +
 433014094785046118400 a^10 b^10 + 27063380924065382400 a^12 b^10 +
 966549318716620800 a^14 b^10 + 15102333104947200 a^16 b^10 +
 15835944037853115187200 b^12 + 15835944037853115187200 a^2 b^12 +
 6928225516560737894400 a^4 b^12 + 1732056379140184473600 a^6 b^12 +
 270633809240653824000 a^8 b^12 + 27063380924065382400 a^10 b^12 +
 1691461307754086400 a^12 b^12 + 60409332419788800 a^14 b^12 +
 943895819059200 a^16 b^12 + 565569429923325542400 b^14 +
 565569429923325542400 a^2 b^14 + 247436625591454924800 a^4 b^14 +
 61859156397863731200 a^6 b^14 + 9665493187166208000 a^8 b^14 +
 966549318716620800 a^10 b^14 + 60409332419788800 a^12 b^14 +
 2157476157849600 a^14 b^14 + 33710564966400 a^16 b^14 +
 8837022342551961600 b^16 + 8837022342551961600 a^2 b^16 +
 3866197274866483200 a^4 b^16 + 966549318716620800 a^6 b^16 +
 151023331049472000 a^8 b^16 + 15102333104947200 a^10 b^16 +
 943895819059200 a^12 b^16 + 33710564966400 a^14 b^16 +
 526727577600 a^16 b^16 + 1178507184834162676727808000 c^2 +
 1189545780444262836889190400 a^2 c^2 +
 525249054431071582342348800 a^4 c^2 +
 132516304956266486169600000 a^6 c^2 +
 20893545903372971802624000 a^8 c^2 +
 2108116095058273370112000 a^10 c^2 +
 132928236194037025996800 a^12 c^2 +
 4789200094813067673600 a^14 c^2 + 75482899175964672000 a^16 c^2 +
 1189545780444262836889190400 b^2 c^2 +
 1200259315024514565695078400 a^2 b^2 c^2 +
 529794886474584937621094400 a^4 b^2 c^2 +
 133617646849611109131878400 a^6 b^2 c^2 +
 21060177855265766375424000 a^8 b^2 c^2 +
 2124237435795556584652800 a^10 b^2 c^2 +
 133902167677007088844800 a^12 b^2 c^2 +
 4822789121955908812800 a^14 b^2 c^2 +
 75989186914340044800 a^16 b^2 c^2 +
 525249054431071582342348800 b^4 c^2 +
 529794886474584937621094400 a^2 b^4 c^2 +
 233771931867713966073446400 a^4 b^4 c^2 +
 58939117104710746216857600 a^6 b^4 c^2 +
 9286674706485076819968000 a^8 b^4 c^2 +
 936402906520464497049600 a^10 b^4 c^2 +
 59008125773560990924800 a^12 b^4 c^2 +
 2124662473700986060800 a^14 b^4 c^2 +
 33466770160562995200 a^16 b^4 c^2 +
 132516304956266486169600000 b^6 c^2 +
 133617646849611109131878400 a^2 b^6 c^2 +
 58939117104710746216857600 a^4 b^6 c^2 +
 14855018310596061403545600 a^6 b^6 c^2 +
 2339866754286146813952000 a^8 b^6 c^2 +
 235861969166930018304000 a^10 b^6 c^2 +
 14858471394875591884800 a^12 b^6 c^2 +
 534837935004136243200 a^14 b^6 c^2 +
 8422067761525555200 a^16 b^6 c^2 +
 20893545903372971802624000 b^8 c^2 +
 21060177855265766375424000 a^2 b^8 c^2 +
 9286674706485076819968000 a^4 b^8 c^2 +
 2339866754286146813952000 a^6 b^8 c^2 +
 368445497799113441280000 a^8 b^8 c^2 +
 37128468304778231808000 a^10 b^8 c^2 +
 2338260850645106688000 a^12 b^8 c^2 +
 84142110929780736000 a^14 b^8 c^2 +
 1324600466079744000 a^16 b^8 c^2 +
 2108116095058273370112000 b^10 c^2 +
 2124237435795556584652800 a^2 b^10 c^2 +
 936402906520464497049600 a^4 b^10 c^2 +
 235861969166930018304000 a^6 b^10 c^2 +
 37128468304778231808000 a^8 b^10 c^2 +
 3740334540329189376000 a^10 b^10 c^2 +
 235487899861711257600 a^12 b^10 c^2 +
 8471567863509811200 a^14 b^10 c^2 +
 133325284442112000 a^16 b^10 c^2 +
 132928236194037025996800 b^12 c^2 +
 133902167677007088844800 a^2 b^12 c^2 +
 59008125773560990924800 a^4 b^12 c^2 +
 14858471394875591884800 a^6 b^12 c^2 +
 2338260850645106688000 a^8 b^12 c^2 +
 235487899861711257600 a^10 b^12 c^2 +
 14821816246033305600 a^12 b^12 c^2 +
 533057065377177600 a^14 b^12 c^2 + 8386907642265600 a^16 b^12 c^2 +
 4789200094813067673600 b^14 c^2 +
 4822789121955908812800 a^2 b^14 c^2 +
 2124662473700986060800 a^4 b^14 c^2 +
 534837935004136243200 a^6 b^14 c^2 +
 84142110929780736000 a^8 b^14 c^2 +
 8471567863509811200 a^10 b^14 c^2 +
 533057065377177600 a^12 b^14 c^2 + 19165729108377600 a^14 b^14 c^2 +
 301463750246400 a^16 b^14 c^2 + 75482899175964672000 b^16 c^2 +
 75989186914340044800 a^2 b^16 c^2 +
 33466770160562995200 a^4 b^16 c^2 +
 8422067761525555200 a^6 b^16 c^2 +
 1324600466079744000 a^8 b^16 c^2 +
 133325284442112000 a^10 b^16 c^2 + 8386907642265600 a^12 b^16 c^2 +
 301463750246400 a^14 b^16 c^2 + 4740548198400 a^16 b^16 c^2 +
 4205418585929895618690416640 c^4 +
 4283789227629131992808816640 a^2 c^4 +
 1908717408260566013957898240 a^4 c^4 +
 485887373715697908622295040 a^6 c^4 +
 77291173706813363611238400 a^8 c^4 +
 7867308845138505608724480 a^10 c^4 +
 500410289854955472814080 a^12 c^4 +
 18185010598892147834880 a^14 c^4 + 289071888149124218880 a^16 c^4 +
 4283789227629131992808816640 b^2 c^4 +
 4359668495142875067894988800 a^2 b^2 c^4 +
 1940827262109243233733181440 a^4 b^2 c^4 +
 493643668385112855360307200 a^6 b^2 c^4 +
 78460828160884078333132800 a^8 b^2 c^4 +
 7980057985327289786695680 a^10 b^2 c^4 +
 507194239139627139072000 a^12 b^2 c^4 +
 18417929327697859706880 a^14 b^2 c^4 +
 292565273543914291200 a^16 b^2 c^4 +
 1908717408260566013957898240 b^4 c^4 +
 1940827262109243233733181440 a^2 b^4 c^4 +
 863281484732864376172707840 a^4 b^4 c^4 +
 219393134698552253837475840 a^6 b^4 c^4 +
 34843167553001087198822400 a^8 b^4 c^4 +
 3541090376643440730439680 a^10 b^4 c^4 +
 224896162105929687367680 a^12 b^4 c^4 +
 8160851353692882862080 a^14 b^4 c^4 +
 129543066926689812480 a^16 b^4 c^4 +
 485887373715697908622295040 b^6 c^4 +
 493643668385112855360307200 a^2 b^6 c^4 +
 219393134698552253837475840 a^4 b^6 c^4 +
 55711953835203493468569600 a^6 b^6 c^4 +
 8841138548643438841036800 a^8 b^6 c^4 +
 897848160984087637524480 a^10 b^6 c^4 +
 56981387314642865356800 a^12 b^6 c^4 +
 2066240988872630599680 a^14 b^6 c^4 +
 32776557579730944000 a^16 b^6 c^4 +
 77291173706813363611238400 b^8 c^4 +
 78460828160884078333132800 a^2 b^8 c^4 +
 34843167553001087198822400 a^4 b^8 c^4 +
 8841138548643438841036800 a^6 b^8 c^4 +
 1401983290824405811200000 a^8 b^8 c^4 +
 142272741033571044556800 a^10 b^8 c^4 +
 9022876604746570137600 a^12 b^8 c^4 +
 326960229932059852800 a^14 b^8 c^4 +
 5183078021568921600 a^16 b^8 c^4 +
 7867308845138505608724480 b^10 c^4 +
 7980057985327289786695680 a^2 b^10 c^4 +
 3541090376643440730439680 a^4 b^10 c^4 +
 897848160984087637524480 a^6 b^10 c^4 +
 142272741033571044556800 a^8 b^10 c^4 +
 14427574047705645711360 a^10 b^10 c^4 +
 914361084276374568960 a^12 b^10 c^4 +
 33111348226576220160 a^14 b^10 c^4 +
 524549875103170560 a^16 b^10 c^4 +
 500410289854955472814080 b^12 c^4 +
 507194239139627139072000 a^2 b^12 c^4 +
 224896162105929687367680 a^4 b^12 c^4 +
 56981387314642865356800 a^6 b^12 c^4 +
 9022876604746570137600 a^8 b^12 c^4 +
 914361084276374568960 a^10 b^12 c^4 +
 57909653012277964800 a^12 b^12 c^4 +
 2095686946952847360 a^14 b^12 c^4 +
 33178745688883200 a^16 b^12 c^4 + 18185010598892147834880 b^14 c^4 +
 18417929327697859706880 a^2 b^14 c^4 +
 8160851353692882862080 a^4 b^14 c^4 +
 2066240988872630599680 a^6 b^14 c^4 +
 326960229932059852800 a^8 b^14 c^4 +
 33111348226576220160 a^10 b^14 c^4 +
 2095686946952847360 a^12 b^14 c^4 +
 75792303051816960 a^14 b^14 c^4 + 1199191897128960 a^16 b^14 c^4 +
 289071888149124218880 b^16 c^4 +
 292565273543914291200 a^2 b^16 c^4 +
 129543066926689812480 a^4 b^16 c^4 +
 32776557579730944000 a^6 b^16 c^4 +
 5183078021568921600 a^8 b^16 c^4 +
 524549875103170560 a^10 b^16 c^4 + 33178745688883200 a^12 b^16 c^4 +
 1199191897128960 a^14 b^16 c^4 + 18962192793600 a^16 b^16 c^4 +
 9099687640450155263180144640 c^6 +
 9328545438479094763770347520 a^2 c^6 +
 4183431927861844189359636480 a^4 c^6 +
 1071928757187759243201085440 a^6 c^6 +
 171644911034162330625638400 a^8 c^6 +
 17588428904315521943470080 a^10 c^6 +
 1126299084549332077117440 a^12 c^6 +
 41209030285424811048960 a^14 c^6 + 659566638832004628480 a^16 c^6 +
 9328545438479094763770347520 b^2 c^6 +
 9548456362608757308433367040 a^2 b^2 c^6 +
 4275723411923189828325212160 a^4 b^2 c^6 +
 1094020420323366888705884160 a^6 b^2 c^6 +
 174943213236931741404364800 a^8 b^2 c^6 +
 17902886471689191334871040 a^10 b^2 c^6 +
 1144990602311449423380480 a^12 b^2 c^6 +
 41842178448700232171520 a^14 b^2 c^6 +
 668921265389064683520 a^16 b^2 c^6 +
 4183431927861844189359636480 b^4 c^6 +
 4275723411923189828325212160 a^2 b^4 c^6 +
 1911906715172373635530752000 a^4 b^4 c^6 +
 488523541416471233581547520 a^6 b^4 c^6 +
 78015659142786066048614400 a^8 b^4 c^6 +
 7973579571090166209576960 a^10 b^4 c^6 +
 509328658192274549637120 a^12 b^4 c^6 +
 18590714330907672576000 a^14 b^4 c^6 +
 296866673570390999040 a^16 b^4 c^6 +
 1071928757187759243201085440 b^6 c^6 +
 1094020420323366888705884160 a^2 b^6 c^6 +
 488523541416471233581547520 a^4 b^6 c^6 +
 124660008517518226567987200 a^6 b^6 c^6 +
 19882257503035609763020800 a^8 b^6 c^6 +
 2029543733065706728980480 a^10 b^6 c^6 +
 129485992445516216401920 a^12 b^6 c^6 +
 4720829210372364042240 a^14 b^6 c^6 +
 75300636585413836800 a^16 b^6 c^6 +
 171644911034162330625638400 b^8 c^6 +
 174943213236931741404364800 a^2 b^8 c^6 +
 78015659142786066048614400 a^4 b^8 c^6 +
 19882257503035609763020800 a^6 b^8 c^6 +
 3167114835118667489280000 a^8 b^8 c^6 +
 322903749309701343436800 a^10 b^8 c^6 +
 20577407598801401241600 a^12 b^8 c^6 +
 749368863900931276800 a^14 b^8 c^6 +
 11939920073677209600 a^16 b^8 c^6 +
 17588428904315521943470080 b^10 c^6 +
 17902886471689191334871040 a^2 b^10 c^6 +
 7973579571090166209576960 a^4 b^10 c^6 +
 2029543733065706728980480 a^6 b^10 c^6 +
 322903749309701343436800 a^8 b^10 c^6 +
 32883202299738670709760 a^10 b^10 c^6 +
 2093140483135599098880 a^12 b^10 c^6 +
 76142089323605329920 a^14 b^10 c^6 +
 1211899133470556160 a^16 b^10 c^6 +
 1126299084549332077117440 b^12 c^6 +
 1144990602311449423380480 a^2 b^12 c^6 +
 509328658192274549637120 a^4 b^12 c^6 +
 129485992445516216401920 a^6 b^12 c^6 +
 20577407598801401241600 a^8 b^12 c^6 +
 2093140483135599098880 a^10 b^12 c^6 +
 133089176877671669760 a^12 b^12 c^6 +
 4836180421689200640 a^14 b^12 c^6 +
 76893677039370240 a^16 b^12 c^6 + 41209030285424811048960 b^14 c^6 +
 41842178448700232171520 a^2 b^14 c^6 +
 18590714330907672576000 a^4 b^14 c^6 +
 4720829210372364042240 a^6 b^14 c^6 +
 749368863900931276800 a^8 b^14 c^6 +
 76142089323605329920 a^10 b^14 c^6 +
 4836180421689200640 a^12 b^14 c^6 +
 175552679006208000 a^14 b^14 c^6 + 2788384179732480 a^16 b^14 c^6 +
 659566638832004628480 b^16 c^6 +
 668921265389064683520 a^2 b^16 c^6 +
 296866673570390999040 a^4 b^16 c^6 +
 75300636585413836800 a^6 b^16 c^6 +
 11939920073677209600 a^8 b^16 c^6 +
 1211899133470556160 a^10 b^16 c^6 +
 76893677039370240 a^12 b^16 c^6 + 2788384179732480 a^14 b^16 c^6 +
 44245116518400 a^16 b^16 c^6 + 13763892782333028209308729344 c^8 +
 14103109556141872838023839744 a^2 c^8 +
 6324673191404469319070908416 a^4 c^8 +
 1621383551616502173007872000 a^6 c^8 +
 259876470921531193389219840 a^8 c^8 +
 26666910412308406377381888 a^10 c^8 +
 1710780105051900475342848 a^12 c^8 +
 62734240818043847442432 a^14 c^8 + 1006726746667147591680 a^16 c^8 +
 14103109556141872838023839744 b^2 c^8 +
 14422290827620942671903719424 a^2 b^2 c^8 +
 6455495906337291520375259136 a^4 b^2 c^8 +
 1651871032509760639504220160 a^6 b^2 c^8 +
 264291445364500788632616960 a^8 b^2 c^8 +
 27073328183710343126581248 a^10 b^2 c^8 +
 1733975602447550710284288 a^12 b^2 c^8 +
 63483411143987520602112 a^14 b^2 c^8 +
 1017187102704449617920 a^16 b^2 c^8 +
 6324673191404469319070908416 b^4 c^8 +
 6455495906337291520375259136 a^2 b^4 c^8 +
 2884161073564195292272656384 a^4 b^4 c^8 +
 736687926794666290603622400 a^6 b^4 c^8 +
 117660716361730324841103360 a^8 b^4 c^8 +
 12032471032357082048888832 a^10 b^4 c^8 +
 769387376904868668506112 a^12 b^4 c^8 +
 28123897552066588704768 a^14 b^4 c^8 +
 449939007070399365120 a^16 b^4 c^8 +
 1621383551616502173007872000 b^6 c^8 +
 1651871032509760639504220160 a^2 b^6 c^8 +
 736687926794666290603622400 a^4 b^6 c^8 +
 187838478195020353744404480 a^6 b^6 c^8 +
 29949502278670934707077120 a^8 b^6 c^8 +
 3057671297003307423989760 a^10 b^6 c^8 +
 195200133636469041561600 a^12 b^6 c^8 +
 7124120769221169315840 a^14 b^6 c^8 +
 113802566745088327680 a^16 b^6 c^8 +
 259876470921531193389219840 b^8 c^8 +
 264291445364500788632616960 a^2 b^8 c^8 +
 117660716361730324841103360 a^4 b^8 c^8 +
 29949502278670934707077120 a^6 b^8 c^8 +
 4767259308474131471923200 a^8 b^8 c^8 +
 485916731439371525422080 a^10 b^8 c^8 +
 30971447743694520222720 a^12 b^8 c^8 +
 1128605181568794961920 a^14 b^8 c^8 +
 18001625313332183040 a^16 b^8 c^8 +
 26666910412308406377381888 b^10 c^8 +
 27073328183710343126581248 a^2 b^10 c^8 +
 12032471032357082048888832 a^4 b^10 c^8 +
 3057671297003307423989760 a^6 b^10 c^8 +
 485916731439371525422080 a^8 b^10 c^8 +
 49449484262486839781376 a^10 b^10 c^8 +
 3146907814019581602816 a^12 b^10 c^8 +
 114499517472837310464 a^14 b^10 c^8 +
 1823600021515130880 a^16 b^10 c^8 +
 1710780105051900475342848 b^12 c^8 +
 1733975602447550710284288 a^2 b^12 c^8 +
 769387376904868668506112 a^4 b^12 c^8 +
 195200133636469041561600 a^6 b^12 c^8 +
 30971447743694520222720 a^8 b^12 c^8 +
 3146907814019581602816 a^10 b^12 c^8 +
 199959513321408454656 a^12 b^12 c^8 +
 7264597477448392704 a^14 b^12 c^8 +
 115531967580825600 a^16 b^12 c^8 +
 62734240818043847442432 b^14 c^8 +
 63483411143987520602112 a^2 b^14 c^8 +
 28123897552066588704768 a^4 b^14 c^8 +
 7124120769221169315840 a^6 b^14 c^8 +
 1128605181568794961920 a^8 b^14 c^8 +
 114499517472837310464 a^10 b^14 c^8 +
 7264597477448392704 a^12 b^14 c^8 +
 263537599896539136 a^14 b^14 c^8 + 4185114339409920 a^16 b^14 c^8 +
 1006726746667147591680 b^16 c^8 +
 1017187102704449617920 a^2 b^16 c^8 +
 449939007070399365120 a^4 b^16 c^8 +
 113802566745088327680 a^6 b^16 c^8 +
 18001625313332183040 a^8 b^16 c^8 +
 1823600021515130880 a^10 b^16 c^8 +
 115531967580825600 a^12 b^16 c^8 + 4185114339409920 a^14 b^16 c^8 +
 66367674777600 a^16 b^16 c^8 + 15756690096587353153351974912 c^10 +
 15989729995385305254359728128 a^2 c^10 +
 7107051537369899368354676736 a^4 c^10 +
 1807142753962795067114913792 a^6 c^10 +
 287517911047290105220104192 a^8 c^10 +
 29309044796012575693209600 a^10 c^10 +
 1869372934965939939508224 a^12 c^10 +
 68205782668786658574336 a^14 c^10 +
 1089891946405907398656 a^16 c^10 +
 15989729995385305254359728128 b^2 c^10 +
 16192078017544193251671539712 a^2 b^2 c^10 +
 7181976429540331818063888384 a^4 b^2 c^10 +
 1822429215693235188274249728 a^6 b^2 c^10 +
 289361208254518886947872768 a^8 b^2 c^10 +
 29438150900192895538298880 a^10 b^2 c^10 +
 1873946422946120358076416 a^12 b^2 c^10 +
 68242692896110318731264 a^14 b^2 c^10 +
 1088468406916298440704 a^16 b^2 c^10 +
 7107051537369899368354676736 b^4 c^10 +
 7181976429540331818063888384 a^2 b^4 c^10 +
 3178937462854661344852033536 a^4 b^4 c^10 +
 804992498756807975455500288 a^6 b^4 c^10 +
 127554011929964148274658304 a^8 b^4 c^10 +
 12950585392186199663345664 a^10 b^4 c^10 +
 822766020254035072717824 a^12 b^4 c^10 +
 29904297569975913154560 a^14 b^4 c^10 +
 476071402416271441920 a^16 b^4 c^10 +
 1807142753962795067114913792 b^6 c^10 +
 1822429215693235188274249728 a^2 b^6 c^10 +
 804992498756807975455500288 a^4 b^6 c^10 +
 203426372326458213545598720 a^6 b^6 c^10 +
 32167738934047077133175424 a^8 b^6 c^10 +
 3259381783834845650128896 a^10 b^6 c^10 +
 206658376561719306306432 a^12 b^6 c^10 +
 7496446165027505167104 a^14 b^6 c^10 +
 119112187839853307904 a^16 b^6 c^10 +
 287517911047290105220104192 b^8 c^10 +
 289361208254518886947872768 a^2 b^8 c^10 +
 127554011929964148274658304 a^4 b^8 c^10 +
 32167738934047077133175424 a^6 b^8 c^10 +
 5076286611527597478163200 a^8 b^8 c^10 +
 513307396459161216428544 a^10 b^8 c^10 +
 32480254759253125136256 a^12 b^8 c^10 +
 1175862257498290053888 a^14 b^8 c^10 +
 18646843045199943168 a^16 b^8 c^10 +
 29309044796012575693209600 b^10 c^10 +
 29438150900192895538298880 a^2 b^10 c^10 +
 12950585392186199663345664 a^4 b^10 c^10 +
 3259381783834845650128896 a^6 b^10 c^10 +
 513307396459161216428544 a^8 b^10 c^10 +
 51799537371501017102592 a^10 b^10 c^10 +
 3271057129163369419776 a^12 b^10 c^10 +
 118182179400091301376 a^14 b^10 c^10 +
 1870412179680688896 a^16 b^10 c^10 +
 1869372934965939939508224 b^12 c^10 +
 1873946422946120358076416 a^2 b^12 c^10 +
 822766020254035072717824 a^4 b^12 c^10 +
 206658376561719306306432 a^6 b^12 c^10 +
 32480254759253125136256 a^8 b^12 c^10 +
 3271057129163369419776 a^10 b^12 c^10 +
 206143663367851808256 a^12 b^12 c^10 +
 7432858707630370560 a^14 b^12 c^10 +
 117400198041907200 a^16 b^12 c^10 +
 68205782668786658574336 b^14 c^10 +
 68242692896110318731264 a^2 b^14 c^10 +
 29904297569975913154560 a^4 b^14 c^10 +
 7496446165027505167104 a^6 b^14 c^10 +
 1175862257498290053888 a^8 b^14 c^10 +
 118182179400091301376 a^10 b^14 c^10 +
 7432858707630370560 a^12 b^14 c^10 +
 267461983098971136 a^14 b^14 c^10 +
 4215962893476096 a^16 b^14 c^10 + 1089891946405907398656 b^16 c^10 +
 1088468406916298440704 a^2 b^16 c^10 +
 476071402416271441920 a^4 b^16 c^10 +
 119112187839853307904 a^6 b^16 c^10 +
 18646843045199943168 a^8 b^16 c^10 +
 1870412179680688896 a^10 b^16 c^10 +
 117400198041907200 a^12 b^16 c^10 +
 4215962893476096 a^14 b^16 c^10 + 66320921028096 a^16 b^16 c^10 +
 13742780160706977462591946752 c^12 +
 13734790680534247390247387136 a^2 c^12 +
 6013639640497792371342704640 a^4 c^12 +
 1506752364669814563319775232 a^6 c^12 +
 236313937135391341691437056 a^8 c^12 +
 23758195326683836579577856 a^10 c^12 +
 1495358873048699714764800 a^12 c^12 +
 53875998037836790824960 a^14 c^12 +
 850753640021892268032 a^16 c^12 +
 13734790680534247390247387136 b^2 c^12 +
 13694237638249619259888500736 a^2 b^2 c^12 +
 5981317532153586870237462528 a^4 b^2 c^12 +
 1494941162577765259899838464 a^6 b^2 c^12 +
 233871812075639679059865600 a^8 b^2 c^12 +
 23452919693787134374871040 a^10 b^2 c^12 +
 1472373465380941577932800 a^12 b^2 c^12 +
 52912361496480338829312 a^14 b^2 c^12 +
 833415323489358446592 a^16 b^2 c^12 +
 6013639640497792371342704640 b^4 c^12 +
 5981317532153586870237462528 a^2 b^4 c^12 +
 2605962317595385182145634304 a^4 b^4 c^12 +
 649650422392580175465194496 a^6 b^4 c^12 +
 101366414409614879861761536 a^8 b^4 c^12 +
 10138058018991001950179328 a^10 b^4 c^12 +
 634751999643330763061760 a^12 b^4 c^12 +
 22749099031565413493760 a^14 b^4 c^12 +
 357345976703529074688 a^16 b^4 c^12 +
 1506752364669814563319775232 b^6 c^12 +
 1494941162577765259899838464 a^2 b^6 c^12 +
 649650422392580175465194496 a^4 b^6 c^12 +
 161525152356629937695561088 a^6 b^6 c^12 +
 25134573894694186953665088 a^8 b^6 c^12 +
 2506812633892787236867584 a^10 b^6 c^12 +
 156508807363778253191616 a^12 b^6 c^12 +
 5593073416898433759360 a^14 b^6 c^12 +
 87602730332509550592 a^16 b^6 c^12 +
 236313937135391341691437056 b^8 c^12 +
 233871812075639679059865600 a^2 b^8 c^12 +
 101366414409614879861761536 a^4 b^8 c^12 +
 25134573894694186953665088 a^6 b^8 c^12 +
 3900132403002524747212128 a^8 b^8 c^12 +
 387855110346978839463936 a^10 b^8 c^12 +
 24143179046172411717792 a^12 b^8 c^12 +
 860177785265876288448 a^14 b^8 c^12 +
 13431301077326598144 a^16 b^8 c^12 +
 23758195326683836579577856 b^10 c^12 +
 23452919693787134374871040 a^2 b^10 c^12 +
 10138058018991001950179328 a^4 b^10 c^12 +
 2506812633892787236867584 a^6 b^10 c^12 +
 387855110346978839463936 a^8 b^10 c^12 +
 38454996493316833239168 a^10 b^10 c^12 +
 2386317482377876081152 a^12 b^10 c^12 +
 84749111716850193024 a^14 b^10 c^12 +
 1319007260073823488 a^16 b^10 c^12 +
 1495358873048699714764800 b^12 c^12 +
 1472373465380941577932800 a^2 b^12 c^12 +
 634751999643330763061760 a^4 b^12 c^12 +
 156508807363778253191616 a^6 b^12 c^12 +
 24143179046172411717792 a^8 b^12 c^12 +
 2386317482377876081152 a^10 b^12 c^12 +
 147604573564618136544 a^12 b^12 c^12 +
 5224598076192918336 a^14 b^12 c^12 +
 81033891718211712 a^16 b^12 c^12 +
 53875998037836790824960 b^14 c^12 +
 52912361496480338829312 a^2 b^14 c^12 +
 22749099031565413493760 a^4 b^14 c^12 +
 5593073416898433759360 a^6 b^14 c^12 +
 860177785265876288448 a^8 b^14 c^12 +
 84749111716850193024 a^10 b^14 c^12 +
 5224598076192918336 a^12 b^14 c^12 +
 184284235532203776 a^14 b^14 c^12 +
 2847900669804288 a^16 b^14 c^12 + 850753640021892268032 b^16 c^12 +
 833415323489358446592 a^2 b^16 c^12 +
 357345976703529074688 a^4 b^16 c^12 +
 87602730332509550592 a^6 b^16 c^12 +
 13431301077326598144 a^8 b^16 c^12 +
 1319007260073823488 a^10 b^16 c^12 +
 81033891718211712 a^12 b^16 c^12 + 2847900669804288 a^14 b^16 c^12 +
 43843813501824 a^16 b^16 c^12 + 8484085605964595841947271168 c^14 +
 8357492172097111418258817024 a^2 c^14 +
 3604367987313881249192306688 a^4 c^14 +
 889019442515791385838240768 a^6 c^14 +
 137185534904462043701732352 a^8 c^14 +
 13564160503520378338427904 a^10 c^14 +
 839348202010490600666112 a^12 c^14 +
 29724796948765719951360 a^14 c^14 +
 461342207452546990080 a^16 c^14 +
 8357492172097111418258817024 b^2 c^14 +
 8207865452903159147004665856 a^2 b^2 c^14 +
 3528439047268155462981385728 a^4 b^2 c^14 +
 867318322578665412093278976 a^6 b^2 c^14 +
 133353272377956693875952384 a^8 b^2 c^14 +
 13135057959056583273179904 a^10 b^2 c^14 +
 809549935322361155516160 a^12 b^2 c^14 +
 28550030654631813917184 a^14 b^2 c^14 +
 441192672894941614080 a^16 b^2 c^14 +
 3604367987313881249192306688 b^4 c^14 +
 3528439047268155462981385728 a^2 b^4 c^14 +
 1511588064954239822516718240 a^4 b^4 c^14 +
 370188705526510989452094864 a^6 b^4 c^14 +
 56693693484796908192803856 a^8 b^4 c^14 +
 5560845258313469790409104 a^10 b^4 c^14 +
 341210205056727472731984 a^12 b^4 c^14 +
 11977026851892969939168 a^14 b^4 c^14 +
 184176695116255380480 a^16 b^4 c^14 +
 889019442515791385838240768 b^6 c^14 +
 867318322578665412093278976 a^2 b^6 c^14 +
 370188705526510989452094864 a^4 b^6 c^14 +
 90298321408963532583380520 a^6 b^6 c^14 +
 13769710164930643317681480 a^8 b^6 c^14 +
 1344392366280131844404712 a^10 b^6 c^14 +
 82084398036284935679400 a^12 b^6 c^14 +
 2866143529903402417968 a^14 b^6 c^14 +
 43828068120342762816 a^16 b^6 c^14 +
 137185534904462043701732352 b^8 c^14 +
 133353272377956693875952384 a^2 b^8 c^14 +
 56693693484796908192803856 a^4 b^8 c^14 +
 13769710164930643317681480 a^6 b^8 c^14 +
 2089967517339107717532312 a^8 b^8 c^14 +
 203018703620992356544488 a^10 b^8 c^14 +
 12327705559305018120024 a^12 b^8 c^14 +
 427897719198467689392 a^14 b^8 c^14 +
 6501573696280797216 a^16 b^8 c^14 +
 13564160503520378338427904 b^10 c^14 +
 13135057959056583273179904 a^2 b^10 c^14 +
 5560845258313469790409104 a^4 b^10 c^14 +
 1344392366280131844404712 a^6 b^10 c^14 +
 203018703620992356544488 a^8 b^10 c^14 +
 19611469461279780410472 a^10 b^10 c^14 +
 1183582106217037189800 a^12 b^10 c^14 +
 40808010742778049264 a^14 b^10 c^14 +
 615523644745847712 a^16 b^10 c^14 +
 839348202010490600666112 b^12 c^14 +
 809549935322361155516160 a^2 b^12 c^14 +
 341210205056727472731984 a^4 b^12 c^14 +
 82084398036284935679400 a^6 b^12 c^14 +
 12327705559305018120024 a^8 b^12 c^14 +
 1183582106217037189800 a^10 b^12 c^14 +
 70946214174613099224 a^12 b^12 c^14 +
 2427659157847862160 a^14 b^12 c^14 +
 36310692234720672 a^16 b^12 c^14 +
 29724796948765719951360 b^14 c^14 +
 28550030654631813917184 a^2 b^14 c^14 +
 11977026851892969939168 a^4 b^14 c^14 +
 2866143529903402417968 a^6 b^14 c^14 +
 427897719198467689392 a^8 b^14 c^14 +
 40808010742778049264 a^10 b^14 c^14 +
 2427659157847862160 a^12 b^14 c^14 +
 82362003967435680 a^14 b^14 c^14 + 1220012565228576 a^16 b^14 c^14 +
 461342207452546990080 b^16 c^14 +
 441192672894941614080 a^2 b^16 c^14 +
 184176695116255380480 a^4 b^16 c^14 +
 43828068120342762816 a^6 b^16 c^14 +
 6501573696280797216 a^8 b^16 c^14 +
 615523644745847712 a^10 b^16 c^14 +
 36310692234720672 a^12 b^16 c^14 + 1220012565228576 a^14 b^16 c^14 +
 17870437083456 a^16 b^16 c^14 + 3200751090373242102395424768 c^16 +
 3120375242287438392916818432 a^2 c^16 +
 1330225478940252748792458240 a^4 c^16 +
 323911462321972767356731008 a^6 c^16 +
 49280117828387104634338944 a^8 c^16 +
 4797500490695802288228480 a^10 c^16 +
 291891297134935206092160 a^12 c^16 +
 10149623019693724379904 a^14 c^16 +
 154457830041819429888 a^16 c^16 +
 3120375242287438392916818432 b^2 c^16 +
 3030161292050287147244138304 a^2 b^2 c^16 +
 1286195435060172020126343552 a^4 b^2 c^16 +
 311694922276325132699328720 a^6 b^2 c^16 +
 47170960300917588749225808 a^8 b^2 c^16 +
 4565343041180991743266320 a^10 b^2 c^16 +
 275974042051168432822704 a^12 b^2 c^16 +
 9527886139626428882784 a^14 b^2 c^16 +
 143861832060274976256 a^16 b^2 c^16 +
 1330225478940252748792458240 b^4 c^16 +
 1286195435060172020126343552 a^2 b^4 c^16 +
 543311146926388146093030016 a^4 b^4 c^16 +
 130954385822258762865421376 a^6 b^4 c^16 +
 19698318888574614042787984 a^8 b^4 c^16 +
 1893525464710155261670976 a^10 b^4 c^16 +
 113592215262719748299680 a^12 b^4 c^16 +
 3888262248501440183552 a^14 b^4 c^16 +
 58147768521809800272 a^16 b^4 c^16 +
 323911462321972767356731008 b^6 c^16 +
 311694922276325132699328720 a^2 b^6 c^16 +
 130954385822258762865421376 a^4 b^6 c^16 +
 31370706809371817312198788 a^6 b^6 c^16 +
 4685981621296027144362068 a^8 b^6 c^16 +
 446875008803611741823188 a^10 b^6 c^16 +
 26565366409412000016476 a^12 b^6 c^16 +
 899923145554946102200 a^14 b^6 c^16 +
 13298564543045002416 a^16 b^6 c^16 +
 49280117828387104634338944 b^8 c^16 +
 47170960300917588749225808 a^2 b^8 c^16 +
 19698318888574614042787984 a^4 b^8 c^16 +
 4685981621296027144362068 a^6 b^8 c^16 +
 694340356303195521055816 a^8 b^8 c^16 +
 65598010412755520282948 a^10 b^8 c^16 +
 3857254377684448422964 a^12 b^8 c^16 +
 129006593320829705528 a^14 b^8 c^16 +
 1877879811124852308 a^16 b^8 c^16 +
 4797500490695802288228480 b^10 c^16 +
 4565343041180991743266320 a^2 b^10 c^16 +
 1893525464710155261670976 a^4 b^10 c^16 +
 446875008803611741823188 a^6 b^10 c^16 +
 65598010412755520282948 a^8 b^10 c^16 +
 6128988029114484676228 a^10 b^10 c^16 +
 355646362213282150124 a^12 b^10 c^16 +
 11706081894379726840 a^14 b^10 c^16 +
 167114780727137712 a^16 b^10 c^16 +
 291891297134935206092160 b^12 c^16 +
 275974042051168432822704 a^2 b^12 c^16 +
 113592215262719748299680 a^4 b^12 c^16 +
 26565366409412000016476 a^6 b^12 c^16 +
 3857254377684448422964 a^8 b^12 c^16 +
 355646362213282150124 a^10 b^12 c^16 +
 20303379571582356340 a^12 b^12 c^16 +
 654826493762846312 a^14 b^12 c^16 +
 9109645290581976 a^16 b^12 c^16 +
 10149623019693724379904 b^14 c^16 +
 9527886139626428882784 a^2 b^14 c^16 +
 3888262248501440183552 a^4 b^14 c^16 +
 899923145554946102200 a^6 b^14 c^16 +
 129006593320829705528 a^8 b^14 c^16 +
 11706081894379726840 a^10 b^14 c^16 +
 654826493762846312 a^12 b^14 c^16 +
 20566861605859504 a^14 b^14 c^16 + 276101622567072 a^16 b^14 c^16 +
 154457830041819429888 b^16 c^16 +
 143861832060274976256 a^2 b^16 c^16 +
 58147768521809800272 a^4 b^16 c^16 +
 13298564543045002416 a^6 b^16 c^16 +
 1877879811124852308 a^8 b^16 c^16 +
 167114780727137712 a^10 b^16 c^16 +
 9109645290581976 a^12 b^16 c^16 + 276101622567072 a^14 b^16 c^16 +
 3519912391812 a^16 b^16 c^16 + 542111206611168649565476884 c^18 +
 524982172319000746684642980 a^2 c^18 +
 222001337560468479090155757 a^4 c^18 +
 53537395403205272787323130 a^6 c^18 +
 8052255681877967991413346 a^8 c^18 +
 773366096017811726390520 a^10 c^18 +
 46313544366395388438237 a^12 c^18 +
 1580944718467944344970 a^14 c^18 + 23549293329040239264 a^16 c^18 +
 524982172319000746684642980 b^2 c^18 +
 505918511898221668363740756 a^2 b^2 c^18 +
 212741600552688131718348849 a^4 b^2 c^18 +
 50972731961891258709020250 a^6 b^2 c^18 +
 7609135378345884022816386 a^8 b^2 c^18 +
 724450448170978615578744 a^10 b^2 c^18 +
 42944058987966317470401 a^12 b^2 c^18 +
 1448513766095888967930 a^14 b^2 c^18 +
 21275411352230094024 a^16 b^2 c^18 +
 222001337560468479090155757 b^4 c^18 +
 212741600552688131718348849 a^2 b^4 c^18 +
 355499814394584303038404373/4 a^4 b^4 c^18 +
 42262848142252246126068549/2 a^6 b^4 c^18 +
 6252034377523966127855049/2 a^8 b^4 c^18 +
 294438468687311891683326 a^10 b^4 c^18 +
 68923472314160553238917/4 a^12 b^4 c^18 +
 1144544048340707430381/2 a^14 b^4 c^18 +
 8248665160417831936 a^16 b^4 c^18 +
 53537395403205272787323130 b^6 c^18 +
 50972731961891258709020250 a^2 b^6 c^18 +
 42262848142252246126068549/2 a^4 b^6 c^18 +
 4978460794046178164373355 a^6 b^6 c^18 +
 728384443784834978229167 a^8 b^6 c^18 +
 67690678695978097084444 a^10 b^6 c^18 +
 7792759928012633416213/2 a^12 b^6 c^18 +
 126764727225743994703 a^14 b^6 c^18 +
 1779973287150248186 a^16 b^6 c^18 +
 8052255681877967991413346 b^8 c^18 +
 7609135378345884022816386 a^2 b^8 c^18 +
 6252034377523966127855049/2 a^4 b^8 c^18 +
 728384443784834978229167 a^6 b^8 c^18 +
 105124445925818113419857 a^8 b^8 c^18 +
 9603934732450167217580 a^10 b^8 c^18 +
 1081774023312500243521/2 a^12 b^8 c^18 +
 17102805310330816883 a^14 b^8 c^18 +
 231114800310247544 a^16 b^8 c^18 +
 773366096017811726390520 b^10 c^18 +
 724450448170978615578744 a^2 b^10 c^18 +
 294438468687311891683326 a^4 b^10 c^18 +
 67690678695978097084444 a^6 b^10 c^18 +
 9603934732450167217580 a^8 b^10 c^18 +
 858101874542786873044 a^10 b^10 c^18 +
 46910552481517533290 a^12 b^10 c^18 +
 1423167304273451128 a^14 b^10 c^18 +
 18098762270078228 a^16 b^10 c^18 +
 46313544366395388438237 b^12 c^18 +
 42944058987966317470401 a^2 b^12 c^18 +
 68923472314160553238917/4 a^4 b^12 c^18 +
 7792759928012633416213/2 a^6 b^12 c^18 +
 1081774023312500243521/2 a^8 b^12 c^18 +
 46910552481517533290 a^10 b^12 c^18 +
 9831401448174012965/4 a^12 b^12 c^18 +
 139790538825606949/2 a^14 b^12 c^18 +
 797031849137168 a^16 b^12 c^18 + 1580944718467944344970 b^14 c^18 +
 1448513766095888967930 a^2 b^14 c^18 +
 1144544048340707430381/2 a^4 b^14 c^18 +
 126764727225743994703 a^6 b^14 c^18 +
 17102805310330816883 a^8 b^14 c^18 +
 1423167304273451128 a^10 b^14 c^18 +
 139790538825606949/2 a^12 b^14 c^18 +
 1773890417556559 a^14 b^14 c^18 + 15738263341034 a^16 b^14 c^18 +
 23549293329040239264 b^16 c^18 +
 21275411352230094024 a^2 b^16 c^18 +
 8248665160417831936 a^4 b^16 c^18 +
 1779973287150248186 a^6 b^16 c^18 +
 231114800310247544 a^8 b^16 c^18 +
 18098762270078228 a^10 b^16 c^18 + 797031849137168 a^12 b^16 c^18 +
 15738263341034 a^14 b^16 c^18 + 33866423320 a^16 b^16 c^18
\end{verbatim}

\vskip 0.5cm

\begin{large} \noindent\textbf{Acknowledgements} \end{large} This work was supported by the Natural Sciences and Engineering Research Council of Canada (Nos. 559668-2021 and
4394-2018).

\end{document}